\documentclass[a4paper,11pt]{article}

\input{diagrams}

\newcommand{\cat}[1]{\ensuremath{\mbox{\bfseries {\upshape {#1}}}}}

\newcommand{\cl}[1]{\ensuremath{\mathcal {#1}}}
\newcommand{\bb}[1]{\ensuremath{\mathbb {#1}}}

\newcommand{\lra}{\ensuremath{\longrightarrow}}
\newcommand{\map}[1]{\ensuremath{\stackrel{{#1}}{\lra}}}
\newcommand{\simarr}{\ensuremath{\stackrel{\sim}{\longrightarrow}}}
\newcommand{\sm}{symmetric multicategory{}}
\newcommand{\sms}{symmetric multicategories{}}

\newcommand{\iso}{isomorphism{}}

\newcommand{\qeed}{\hfill $\Box$}

\newtheorem{theorem}{Theorem}[section]
\newtheorem{proposition}[theorem]{Proposition}

\newtheorem{lemma}[theorem]{Lemma}
\newtheorem{definition}[theorem]{Definition}
\newtheorem{definitions}[theorem]{Definitions}

\newtheorem{examples}[theorem]{Examples}
\newtheorem{remarks}[theorem]{Remarks}

\newenvironment{prf}{\vspace{2ex}\begin{sloppypar}{\noindent\upshape
{\bfseries Proof. }}} {{\hspace*{\fill}
$\Box$}\end{sloppypar}\vspace{2ex}}

\newenvironment{prfof}[1]{\vspace{2ex}\begin{sloppypar}{\noindent
\upshape{\bfseries Proof of {#1}. }}} {{\hspace*{\fill}
$\Box$}\end{sloppypar}\vspace{2ex}}

\newcommand{\numroman}{\renewcommand{\labelenumi}{\roman{enumi})}}

\newcommand{\pica}{\begin{center} \input}
\newcommand{\picz}{\end{center}}
\newcommand{\length}[1]{\setlength{\unitlength}{#1}}

\newlength{\leng}
\newlength{\fontleng}
\newcommand{\sunit}{\setlength{\unitlength}{1mm}}
\newcommand{\setleng}[1]{\setlength{\leng}{#1}
    \setlength{\unitlength}{\leng}}
\newcommand{\setunit}[1]{\setlength{\unitlength}{#1}}

\newcommand{\diagc}[3]{

\put(#1,#2){
\setlength{\unitlength}{#3} %

\put(1,2){\line(1,0){5}} %
\put(1,3){\line(1,0){5}} %
\put(1,4){\line(1,0){5}} %
\put(1,1){\line(1,0){5}} %
\put(3.5,-1.5){\line(1,1){4}} %
\put(3.5,6.5){\line(1,-1){4}}}}

\newcommand{\diagf}[3]{

\put(#1,#2){
\setlength{\unitlength}{#3} %

\put(0,0){\line(1,0){4}}  %
\put(0,0){\line(2,3){2}}  %
\put(4,0){\line(-2,3){2}}}}

\newcommand{\diagj}[3]{

\put(#1,#2){
\setlength{\unitlength}{#3} %

\put(0,0){\line(1,0){4}}  %
\qbezier(0,0)(2,2)(4,0) }}

\newcommand{\diagp}[3]{

\put(#1,#2){

\setlength{\unitlength}{#3}  %
\fontsize{#3}{20}

\put(0,0){\makebox(0,0){$\equiv$}}  %
\put(0.35,0){\makebox(0,0){$\rangle$}}}}

\newcommand{\diagr}[3]{

\put(#1,#2){\setlength{\leng}{#3}

\setlength{\unitlength}{1.4\leng}

\qbezier(0.5,0)(2.5,1)(3,4)   %
\put(0.5,-0){\makebox(0,0)[c]{{\tiny $\Delta$}}} }}

\newcommand{\diagad}[3]{

\put(#1,#2){

\setlength{\unitlength}{#3}  %

\qbezier(0.5,0)(2.5,-0.5)(3,-4)   %
\put(3,-4){\makebox(0,0)[t]{{\tiny $\nabla$}}}}}

\newcommand{\diagal}[3]{

\put(#1,#2){
\setlength{\unitlength}{#3} %

\put(0,0){\line(1,0){4}}  %
\put(0,0){\line(2,3){2}}  %
\put(4,0){\line(-2,3){2}} %
\qbezier(0,0)(2,2)(4,0) }}

\newcommand{\diagbj}[3]{

\put(#1,#2){ \setlength{\leng}{#3} %
\setlength{\unitlength}{1.2\leng} %
\setlength{\fontleng}{\leng}   %
\fontsize{\fontleng}{15}

\put(3.5,0){\line(-3,2){3}} %
\put(4,0.5){\line(-3,2){3}} %
\put(3.35,0){\makebox(0,0)[c]{$\underline{\hspace{0.5em}}$}}  %
\put(4.1,0.65){\makebox(0,0)[c]{$\mid$}}  }}

\newcommand{\diagbk}[3]{

\put(#1,#2){ \setlength{\leng}{#3} %
\setlength{\unitlength}{1.2\leng} %
\setlength{\fontleng}{\leng}   %
\fontsize{\fontleng}{15}

\put(-3.5,0){\line(3,2){3}} %
\put(-4,0.5){\line(3,2){3}} %
\put(-3.35,0.0){\makebox(0,0)[c]{$\underline{\hspace{0.5em}}$}}  %
\put(-4,0.65){\makebox(0,0)[c]{$\mid$}}  }}

\sunit

\newcommand{\onetwo}[8]{

\put(#1,#2){ \setlength{\leng}{#3} %
\setlength{\unitlength}{\leng} %

\diagj{0}{0}{6\leng} %
\put(-1,0){\makebox(0,0)[r]{#4}} %
\put(26,0){\makebox(0,0)[l]{#5}} %
\put(12,7){\makebox(0,0)[b]{#6}} %
\put(12,-2){\makebox(0,0)[t]{#7}} %
\put(12,2.5){\makebox(0,0)[c]{$\Downarrow$}}  %
\put(14,2.5){\makebox(0,0)[l]{#8}} }}

\sunit

\newcommand{\onetwob}[8]{

\put(#1,#2){ \setlength{\leng}{#3} %
\setlength{\unitlength}{\leng} %

\diagj{0}{0}{6\leng} %
\put(-1,0){\makebox(0,0)[r]{#4}} %
\put(26,0){\makebox(0,0)[l]{#5}} %
\put(12,7){\makebox(0,0)[b]{#6}} %
\put(12,-2){\makebox(0,0)[t]{#7}} %
\put(12,2.5){\makebox(0,0)[c]{#8}}   }}

\newcommand{\diagbl}[3]{

\put(#1,#2){\setlength{\leng}{#3}
\setlength{\unitlength}{.4\leng} %

\put(0,0){\line(1,2){10}}  %
\put(10,20){\line(1,0){20}}  %
\put(30,20){\line(1,-2){10}}  %
\put(0,0){\line(1,0){40}}}}

\sunit

\newcommand{\nulltwo}[6]{

\put(#1,#2){\setlength{\leng}{#3}
\setlength{\unitlength}{1.2\leng} %
\fontsize{3.5\leng}{15}

\put(5,0){\vector(1,0){10}}  %
\put(10,9){\makebox(0,0)[c]{.}}  %
\put(10,3.5){\makebox(0,0)[c]{$\Downarrow$}}

\fontsize{11pt}{15} %
\put(4,0){\makebox(0,0)[r]{#4}}  %
\put(16.5,0){\makebox(0,0)[l]{#4}}  %
\put(10,-2){\makebox(0,0)[t]{#5}}  %
\put(12,3.5){\makebox(0,0)[l]{#6}}  %

}} %

\sunit

\newcommand{\onecell}[6]{

\put(#1,#2){\setlength{\leng}{#3}
\setlength{\unitlength}{1.2\leng} %

\put(4,0){\makebox(0,0)[r]{#4}}  %
\put(16.5,0){\makebox(0,0)[l]{#5}}  %
\put(10,1.5){\makebox(0,0)[b]{#6}}  %

\put(5,0){\vector(1,0){10}}}}  %

\newcommand{\threetwo}[9]{

\put(0,0){\setlength{\leng}{1.3mm}
\setlength{\unitlength}{\leng} %

\diagbl{0}{0}{\leng}  %

\put(-2,-1){\makebox(0,0)[c]{#1}}  %
\put(2,9){\makebox(0,0)[c]{#2}}  %
\put(15,9){\makebox(0,0)[c]{#3}}  %
\put(19,-1){\makebox(0,0)[c]{#4}}  %

\put(0,5){\makebox(0,0)[c]{#5}}  %
\put(8.5,11){\makebox(0,0)[c]{#6}}  %
\put(17,5){\makebox(0,0)[c]{#7}}  %
\put(8.5,-3){\makebox(0,0)[c]{#8}}  %

\put(8.5,4){\makebox(0,0)[c]{$\Downarrow$}}  %
\put(11,4){\makebox(0,0)[c]{#9}} }} %

\newcommand{\threetwob}[9]{

\put(0,0){\setlength{\leng}{1.3mm}
\setlength{\unitlength}{\leng} %

\diagbl{0}{0}{\leng}  %

\put(-2,-1){\makebox(0,0)[c]{#1}}  %
\put(2,9){\makebox(0,0)[c]{#2}}  %
\put(15,9){\makebox(0,0)[c]{#3}}  %
\put(19,-1){\makebox(0,0)[c]{#4}}  %

\put(0,5){\makebox(0,0)[c]{#5}}  %
\put(8.5,11){\makebox(0,0)[c]{#6}}  %
\put(17,5){\makebox(0,0)[c]{#7}}  %
\put(8.5,-3){\makebox(0,0)[c]{#8}}  %

\put(8.5,4){\makebox(0,0)[c]{#9}}   }} %

\newcommand{\diagbo}[3]{

\put(#1,#2){\setlength{\leng}{#3}
\setlength{\unitlength}{0.4\leng} %

\diagp{68}{12}{6\leng}

\begin{picture}(100,40)
\put(0,0){

\begin{picture}(50,30)      %

\put(0,0){\line(1,2){10}}  %
\put(10,20){\line(1,1){10}}  %
\put(20,30){\line(1,-1){10}}  %
\put(30,20){\line(1,0){15}}  %
\put(45,20){\line(1,-2){5}}  %
\put(40,0){\line(1,1){10}}  %
\put(0,0){\line(1,0){40}}  %
\put(10,20){\line(1,0){20}}  %
\put(30,20){\line(1,-2){10}}  %

\end{picture}}

\put(80,0){
\begin{picture}(50,30)      %

\put(0,0){\line(1,2){10}}  %
\put(10,20){\line(1,1){10}}  %
\put(20,30){\line(1,-1){10}}  %
\put(30,20){\line(1,0){15}}  %
\put(45,20){\line(1,-2){5}}  %
\put(40,0){\line(1,1){10}}  %
\put(0,0){\line(1,0){40}}  %

\end{picture}}

\end{picture}}}

\newcommand{\threethree}[5]{

\put(0,0){\setlength{\leng}{1mm}
\setlength{\unitlength}{0.4\leng} %

\diagp{68}{12}{6\leng} %
\put(70,20){\makebox(0,0)[b]{#5}}

\begin{picture}(100,40)
\put(0,0){

\begin{picture}(50,30)      %

\put(0,0){\line(1,2){10}}  %
\put(10,20){\line(1,1){10}}  %
\put(20,30){\line(1,-1){10}}  %
\put(30,20){\line(1,0){15}}  %
\put(45,20){\line(1,-2){5}}  %
\put(40,0){\line(1,1){10}}  %
\put(0,0){\line(1,0){40}}  %
\put(10,20){\line(1,0){20}}  %
\put(30,20){\line(1,-2){10}}  %
\put(20,10){\makebox(0,0){#3}}  %
\put(42,12){\makebox(0,0){#2}}  %
\put(21,23){\makebox(0,0){#1}}  %

\end{picture}}

\put(80,0){
\begin{picture}(50,30)      %

\put(0,0){\line(1,2){10}}  %
\put(10,20){\line(1,1){10}}  %
\put(20,30){\line(1,-1){10}}  %
\put(30,20){\line(1,0){15}}  %
\put(45,20){\line(1,-2){5}}  %
\put(40,0){\line(1,1){10}}  %
\put(0,0){\line(1,0){40}}  %
\put(25,10){\makebox(0,0){#4}}  %

\end{picture}}

\end{picture}}}

\newcommand{\diagbp}[3]{

\put(#1,#2){\setlength{\leng}{#3}
\setlength{\unitlength}{0.4\leng} %

\begin{picture}(40,30)      %

\put(0,0){\line(1,4){5}}  %
\put(5,20){\line(3,2){15}}  %

\qbezier[2](26,31)(29,28)(32,25) %

\put(35,20){\line(1,-4){5}}  %
\put(0,0){\line(1,0){40}}  %

\put(20,13){\makebox(0,0){$\Downarrow$}}  %

\end{picture}}}

\sunit

\newcommand{\three}[4]{
\put(#1,#2){ \setlength{\leng}{#3} %
\setlength{\unitlength}{\leng} %
\diagp{0}{0}{6\leng}  %
\put(0,4){\makebox(0,0)[b]{#4}} }}

\sunit

\newcommand{\nullthree}[3]{

\put(0,0){ \setlength{\leng}{1mm} %
\setlength{\unitlength}{\leng} %

\onecell{0}{0}{\leng}{}{}{#1}   %
\three{24}{2}{\leng}{#3}  %
\onetwo{31}{0}{0.7\leng}{}{}{#1}{#1}{#2}

 }}

\sunit

\newcommand{\nullthreeb}[3]{

\put(0,0){ \setlength{\leng}{1mm} %
\setlength{\unitlength}{\leng} %

\onecell{0}{0}{\leng}{}{}{#1}   %
\three{24}{2}{\leng}{#3}  %
\onetwob{31}{0}{0.7\leng}{}{}{#1}{#1}{#2}

 }}

\sunit

\newcommand{\twotwo}[7]{

\put(0,0){ \setlength{\leng}{1mm} %
\setlength{\unitlength}{\leng} %
\diagf{0}{0}{4\leng}  %
\put(-1,-1){\makebox(0,0)[tr]{#1}}   
\put(8,14){\makebox(0,0)[b]{#2}}  %
\put(18,-1){\makebox(0,0)[tl]{#3}}  %

\put(3,6){\makebox(0,0)[br]{#4}}   
\put(14,6){\makebox(0,0)[bl]{#5}}  %
\put(8,-2){\makebox(0,0)[t]{#6}}  %

\put(8,5){\makebox(0,0)[c]{$\Downarrow$}}  %
\put(9,5){\makebox(0,0)[l]{#7}}  %
}}

\sunit

\newcommand{\twotwob}[7]{

\put(0,0){ \setlength{\leng}{1mm} %
\setlength{\unitlength}{\leng} %
\diagf{0}{0}{4\leng}  %
\put(-1,-1){\makebox(0,0)[tr]{#1}}   
\put(8,14){\makebox(0,0)[b]{#2}}  %
\put(18,-1){\makebox(0,0)[tl]{#3}}  %

\put(3,6){\makebox(0,0)[br]{#4}}   
\put(14,6){\makebox(0,0)[bl]{#5}}  %
\put(8,-2){\makebox(0,0)[t]{#6}}  %

\put(8,5){\makebox(0,0)[c]{#7}}  %

}}

\newcommand{\assleft}[5]{

\put(0,0){\setlength{\leng}{1.3mm}
\setlength{\unitlength}{\leng} %

\put(0,0){
\setlength{\unitlength}{.4\leng} %

\put(0,0){\line(1,2){10}}  %
\put(10,20){\line(1,0){20}}  %
\put(30,20){\line(1,-2){10}}  %
\put(0,0){\line(1,0){40}}   %
\put(10,20){\line(3,-2){30}}}

\put(0,5){\makebox(0,0)[c]{#1}}  %
\put(8.5,11){\makebox(0,0)[c]{#2}}  %
\put(17,5){\makebox(0,0)[c]{#3}}  %
\put(8.5,-3){\makebox(0,0)[c]{#5}}  %

\put(8.5,3.3){\makebox(0,0)[c]{#4}}   }}

\newcommand{\assright}[5]{

\put(0,0){\setlength{\leng}{1.3mm}
\setlength{\unitlength}{\leng} %

\put(0,0){
\setlength{\unitlength}{.4\leng} %

\put(0,0){\line(1,2){10}}  %
\put(10,20){\line(1,0){20}}  %
\put(30,20){\line(1,-2){10}}  %
\put(0,0){\line(1,0){40}}   %
\put(0,0){\line(3,2){30}}}

\put(0,5){\makebox(0,0)[c]{#1}}  %
\put(8.5,11){\makebox(0,0)[c]{#2}}  %
\put(17,5){\makebox(0,0)[c]{#3}}  %
\put(8.5,-3){\makebox(0,0)[c]{#5}}  %

\put(8.5,3.3){\makebox(0,0)[c]{#4}}   }}

\newcommand{\assrightb}[2]{

\put(0,0){\setlength{\leng}{1.3mm}
\setlength{\unitlength}{\leng} %

\put(0,0){
\setlength{\unitlength}{.4\leng} %

\put(0,0){\line(1,2){10}}  %
\put(10,20){\line(1,0){20}}  %
\put(30,20){\line(1,-2){10}}  %
\put(0,0){\line(1,0){40}}   %
\put(0,0){\line(3,2){30}}}

\put(5,5){\makebox(0,0)[c]{#1}}  %

\put(10,2.5){\makebox(0,0)[c]{#2}}   }}

\newcommand{\assleftb}[2]{

\put(0,0){\setlength{\leng}{1.3mm}
\setlength{\unitlength}{\leng} %

\put(0,0){
\setlength{\unitlength}{.4\leng} %

\put(0,0){\line(1,2){10}}  %
\put(10,20){\line(1,0){20}}  %
\put(30,20){\line(1,-2){10}}  %
\put(0,0){\line(1,0){40}}   %
\put(10,20){\line(3,-2){30}}}

\put(11,5){\makebox(0,0)[c]{#1}}  %

\put(5.5,2.5){\makebox(0,0)[c]{#2}}   }}

\newcommand{\diagbq}[5]{

\put(0,0){\setlength{\leng}{1mm}
\setlength{\unitlength}{\leng} %

\onetwob{0}{0}{\leng}{#1}{#1}{#2}{#3}{#5}
\put(12,21){\makebox(0,0)[c]{.}}  %
\put(12,15.5){\makebox(0,0)[c]{$\Downarrow$}} %
\put(14,15.5){\makebox(0,0)[l]{#4}}  }}

\newcommand{\diagbr}[9]{

\put(0,0){\setlength{\leng}{1mm}
\setlength{\unitlength}{\leng} %

\diagal{0}{0}{4\leng}
\put(-1,-1){\makebox(0,0)[tr]{#1}}   
\put(8,14){\makebox(0,0)[b]{#2}}  %
\put(18,-1){\makebox(0,0)[tl]{#3}}  %

\put(3,6){\makebox(0,0)[br]{#4}}   
\put(14,6){\makebox(0,0)[bl]{#5}}  %
\put(8,-2){\makebox(0,0)[t]{#7}}  %

\put(8,7){\makebox(0,0)[c]{#8}}  %

\put(10,3.3){\makebox(0,0)[bl]{#6}}  %
\put(8,1.6){\makebox(0,0)[c]{#9}}  %
}}

\newcommand{\diagbs}[3]{

\put(#1,#2){
\setlength{\unitlength}{#3} %

\put(0,0){\line(1,0){4}}  %
\put(0,0){\line(2,3){2}}  %
\put(4,0){\line(-2,3){2}} %
\qbezier(2,3)(4.9,2.6)(4,0)}}

\newcommand{\diagbt}[9]{

\put(0,0){\setlength{\leng}{1mm}
\setlength{\unitlength}{\leng} %

\diagbs{0}{0}{4\leng}
\put(-1,-1){\makebox(0,0)[tr]{#1}}   
\put(8,14){\makebox(0,0)[b]{#2}}  %
\put(18,-1){\makebox(0,0)[tl]{#3}}  %

\put(3,6){\makebox(0,0)[br]{#4}}   
\put(12,4){\makebox(0,0)[c]{#6}}  %
\put(8,-2){\makebox(0,0)[t]{#7}}  %

\put(8,6){\makebox(0,0)[c]{#9}}  %

\put(15,6){\makebox(0,0)[c]{#8}}  %
\put(17,9){\makebox(0,0)[bl]{#5}}  %
}}

\newcommand{\diagbu}[3]{

\put(#1,#2){
\setlength{\unitlength}{#3} %

\put(0,0){\line(1,0){4}}  %
\put(0,0){\line(2,3){2}}  %
\put(4,0){\line(-2,3){2}} %
\qbezier(2,3)(-0.9,2.6)(0,0)}}

\newcommand{\diagbv}[9]{

\put(0,0){\setlength{\leng}{1mm}
\setlength{\unitlength}{\leng} %

\diagbu{0}{0}{4\leng}
\put(-1,-1){\makebox(0,0)[tr]{#1}}   
\put(8,14){\makebox(0,0)[b]{#2}}  %
\put(18,-1){\makebox(0,0)[tl]{#3}}  %

\put(3,3){\makebox(0,0)[c]{#5}}   
\put(13,6){\makebox(0,0)[bl]{#6}}  %
\put(8,-2){\makebox(0,0)[t]{#7}}  %

\put(8,6){\makebox(0,0)[c]{#9}}  %

\put(2,7){\makebox(0,0)[c]{#8}}  %
\put(0,8){\makebox(0,0)[br]{#4}}  %
}}

\newcommand{\diagbw}[8]{

\put(0,0){\setlength{\leng}{1mm}
\setlength{\unitlength}{\leng} %

\twotwob{#1}{#2}{#3}{#4}{#5}{#6}{#8} %
\diagbj{-6}{10}{\leng}   %
\put(-9,16){.}   %
\put(-3,13){\makebox(0,0)[bl]{#7}} %
}}

\newcommand{\diagbx}[8]{

\put(0,0){\setlength{\leng}{1mm}
\setlength{\unitlength}{\leng} %

\twotwob{#1}{#2}{#3}{#4}{#5}{#6}{#8} %
\diagbk{22}{10}{\leng}   %
\put(25,15){.}   %
\put(18,13){\makebox(0,0)[br]{#7}} %
}}

\newcommand{\diagby}[2]{

\setlength{\leng}{#1}
\setlength{\unitlength}{\leng} %
\qbezier(0,0)(8,-8)(16,0)      %
\put(8,-2){\makebox(0,0)[c]{#2}} }

\newcommand{\diagbz}[1]{

\put(0,0){\setlength{\leng}{1mm}
\setlength{\unitlength}{\leng} %
\put(0,0){
\setlength{\unitlength}{4\leng} %
\qbezier(0,0)(2,-4)(4,0) }}
\put(8,-6.5){\makebox(0,0)[c]{#1}}  %
}

\newcommand{\diagca}{

\put(0,0){\setlength{\leng}{1mm}
\setlength{\unitlength}{\leng} %
\put(0,0){
\setlength{\unitlength}{4\leng} %
\qbezier(2,3)(6.7,4.1)(4,0)}}}

\newcommand{\diagcb}{

\put(0,0){\setlength{\leng}{1mm}
\setlength{\unitlength}{\leng} %
\put(0,0){
\setlength{\unitlength}{4\leng} %
\qbezier(2,3)(-2.7,4.1)(0,0)}} }

\newcommand{\diagcc}[4]{

\put(0,0){\setlength{\leng}{1.3mm}
\setlength{\unitlength}{\leng} %

\put(0,0){
\setlength{\unitlength}{.4\leng} %
\qbezier(0,0)(20,13)(40,0) %
\put(20,13){\makebox(0,0)[c]{#1}}
\put(20,2.4){\makebox(0,0)[c]{#4}}}  %
\put(10,4){\makebox(0,0)[c]{#3}}  %
\put(5.5,4){\makebox(0,0)[c]{#2}} %
}}

\newcommand{\diagcd}[1]{

\put(0,0){\setlength{\leng}{1.3mm}
\setlength{\unitlength}{\leng} %

\put(0,0){
\setlength{\unitlength}{.4\leng} %
\qbezier(0,0)(20,-13)(40,0) %
\put(20,-3){\makebox(0,0)[c]{#1}}}  %
 }}

\newcommand{\diagce}[1]{

\put(0,0){\setlength{\leng}{1.3mm}
\setlength{\unitlength}{\leng} %

\put(0,0){
\setlength{\unitlength}{.4\leng} %
\qbezier(0,0)(20,-28)(40,0) %
\put(20,-10){\makebox(0,0)[c]{#1}}}  %
 }}

\newcommand{\diagcf}[1]{

\put(0,0){\setlength{\leng}{1.3mm}
\setlength{\unitlength}{\leng} %

\put(0,0){
\setlength{\unitlength}{.4\leng} %
\qbezier(0,0)(-8,18)(10,20) %
\put(2,12){\makebox(0,0)[c]{#1}}}  %
 }}

\newcommand{\diagcg}[1]{

\put(0,0){\setlength{\leng}{1.3mm}
\setlength{\unitlength}{\leng} %

\put(0,0){
\setlength{\unitlength}{.4\leng} %
\qbezier(40,0)(48,18)(30,20) %
\put(38,12){\makebox(0,0)[c]{#1}}}  %
 }}

\newcommand{\diagch}[1]{

\put(0,0){\setlength{\leng}{1.3mm}
\setlength{\unitlength}{\leng} %

\put(0,0){
\setlength{\unitlength}{.4\leng} %
\qbezier(40,0)(58,23)(30,20) %
\put(45,14){\makebox(0,0)[c]{#1}}}  %
 }}

\newcommand{\diagci}[1]{

\put(0,0){\setlength{\leng}{1.3mm}
\setlength{\unitlength}{\leng} %

\put(0,0){
\setlength{\unitlength}{.4\leng} %
\qbezier(10,20)(20,34)(30,20) %
\put(20,23){\makebox(0,0)[c]{#1}}}  %
 }}

\newcommand{\diagcj}[2]{

\put(0,0){\setlength{\leng}{1.3mm}
\setlength{\unitlength}{\leng} %

\put(0,0){
\setlength{\unitlength}{.4\leng} %
\qbezier(0,0)(10,20)(30,20) %
\put(11,16){\makebox(0,0)[c]{{\tiny #1}}} %
\put(15,14){\makebox(0,0)[c]{{\small #2}}} %
}}}

\newcommand{\diagck}[2]{

\put(0,0){\setlength{\leng}{1.3mm}
\setlength{\unitlength}{\leng} %

\put(0,0){
\setlength{\unitlength}{.4\leng} %
\qbezier(40,0)(30,20)(10,20) %
\put(29,16){\makebox(0,0)[c]{{\tiny #1}}} %
\put(25,14){\makebox(0,0)[c]{{\small #2}}} %
}}}

\newcommand{\diagcl}[1]{

\put(0,0){\setlength{\leng}{1.3mm}
\setlength{\unitlength}{\leng} %

\put(0,0){
\setlength{\unitlength}{.4\leng} %
\put(20,28){\makebox(0,0)[c]{$\Downarrow$}} %
\put(20,37){\makebox(0,0)[c]{.}}
\put(25,28){\makebox(0,0)[c]{#1}} %
}}}

\newcommand{\diagcm}{

\put(0,0){\setlength{\leng}{1.3mm}
\setlength{\unitlength}{\leng} %

\put(0,0){
\setlength{\unitlength}{.4\leng} %
\put(10,20){\line(1,1){10}}  %
\put(30,20){\line(-1,1){10}}

 }}}

\newcommand{\mtwo}[5]{

\put(0,0){\setlength{\leng}{1.3mm}
\setlength{\unitlength}{\leng} %

\put(0,0){ \setlength{\unitlength}{.4\leng} %
\put(0,0){\line(1,0){40}}  %
\put(0,0){\line(1,4){5}}  %
\put(5,20){\line(1,1){10}}  %
\put(40,0){\line(-1,4){5}}  %

\put(0,12){\makebox(0,0)[r]{#1}}  %
\put(8,27){\makebox(0,0)[br]{#2}}  %
\put(40,12){\makebox(0,0)[l]{#3}}  %
\put(20,-2){\makebox(0,0)[t]{#4}}  %
\put(20,15){\makebox(0,0)[c]{$\Downarrow$}}  %
\put(23,15){\makebox(0,0)[l]{#5}}  %
\qbezier[2](23,30)(28,27)(33,24)

 }}}

\newcommand{\mtwob}[5]{

\put(0,0){\setlength{\leng}{1.3mm}
\setlength{\unitlength}{\leng} %

\put(0,0){ \setlength{\unitlength}{.4\leng} %
\put(0,0){\line(1,0){40}}  %
\put(0,0){\line(1,4){5}}  %
\put(5,20){\line(1,1){10}}  %
\put(40,0){\line(-1,4){5}}  %

\put(0,12){\makebox(0,0)[r]{#1}}  %
\put(8,27){\makebox(0,0)[br]{#2}}  %
\put(40,12){\makebox(0,0)[l]{#3}}  %
\put(20,-2){\makebox(0,0)[t]{#4}}  %
\put(20,15){\makebox(0,0)[c]{#5}}  %
\qbezier[2](23,30)(28,27)(33,24)

 }}}

\newcommand{\diagco}[2]{

\setlength{\leng}{#1}
\setlength{\unitlength}{\leng} %

\put(0,0){\line(1,0){16}} %
\qbezier(0,0)(8,8)(16,0)  %
\put(8,2){\makebox(0,0)[c]{#2}}}

\newcommand{\onetwoc}[5]{

\put(0,0){ \setlength{\leng}{1mm} %
\setlength{\unitlength}{\leng} %

\diagj{0}{0}{6\leng} %
\put(0,0){ \setlength{\unitlength}{6\leng}
\put(0,0){  %
\qbezier(0,0)(2,-2)(4,0) }}
\put(12,7){\makebox(0,0)[b]{#1}} %
\put(7,0){\makebox(0,0)[b]{#2}} %
\put(12,2.5){\makebox(0,0)[c]{#4}}
\put(12,-7){\makebox(0,0)[t]{#3}} %
\put(12,-3){\makebox(0,0)[c]{#5}}

 }}

\newcommand{\scaleby}[1]{
\setlength{\unitlength}{\leng} %
\qbezier(0,0)(8,-8)(16,0)      %
\put(8,-2){\makebox(0,0)[c]{#1}} }

\newcommand{\scalebv}[9]{

\put(0,0){
\setlength{\unitlength}{\leng} %

\diagbu{0}{0}{4\leng}
\put(-1,-1){\makebox(0,0)[tr]{#1}}   
\put(8,14){\makebox(0,0)[b]{#2}}  %
\put(18,-1){\makebox(0,0)[tl]{#3}}  %

\put(3,3){\makebox(0,0)[c]{#5}}   
\put(13,6){\makebox(0,0)[bl]{#6}}  %
\put(8,-2){\makebox(0,0)[t]{#7}}  %

\put(8,6){\makebox(0,0)[c]{#9}}  %

\put(2,7){\makebox(0,0)[c]{#8}}  %
\put(0,8){\makebox(0,0)[br]{#4}}  %
}}

\newcommand{\scalebt}[9]{

\put(0,0){
\setlength{\unitlength}{\leng} %

\diagbs{0}{0}{4\leng}
\put(-1,-1){\makebox(0,0)[tr]{#1}}   
\put(8,14){\makebox(0,0)[b]{#2}}  %
\put(18,-1){\makebox(0,0)[tl]{#3}}  %

\put(3,6){\makebox(0,0)[br]{#4}}   
\put(13,3.5){\makebox(0,0)[c]{#6}}  %
\put(8,-2){\makebox(0,0)[t]{#7}}  %

\put(8,6){\makebox(0,0)[c]{#9}}  %

\put(15,6){\makebox(0,0)[c]{#8}}  %
\put(17,9){\makebox(0,0)[bl]{#5}}  %
}}

\newcommand{\scalebr}[9]{

\put(0,0){
\setlength{\unitlength}{\leng} %

\diagal{0}{0}{4\leng}
\put(-1,-1){\makebox(0,0)[tr]{#1}}   
\put(8,14){\makebox(0,0)[b]{#2}}  %
\put(18,-1){\makebox(0,0)[tl]{#3}}  %

\put(3,6){\makebox(0,0)[br]{#4}}   
\put(14,6){\makebox(0,0)[bl]{#5}}  %
\put(8,-2){\makebox(0,0)[t]{#7}}  %

\put(8,7){\makebox(0,0)[c]{#8}}  %

\put(10,3){\makebox(0,0)[bl]{#6}}  %
\put(8,1){\makebox(0,0)[c]{#9}}  %
}}

\newcommand{\scaleco}[1]{

\setlength{\unitlength}{\leng} %

\put(0,0){\line(1,0){16}} %
\qbezier(0,0)(8,8)(16,0)  %
\put(8,2){\makebox(0,0)[c]{#1}}}

\newcommand{\scalecp}{

\put(0,0){
\setlength{\unitlength}{\leng} %
\put(0,0){\line(1,0){40}} %
\put(0,0){\line(0,1){20}} %
\put(0,20){\line(2,1){20}} %
\put(20,30){\line(2,-1){20}} %
\put(40,0){\line(0,1){20}}}}

\newcommand{\scalecq}{

\put(0,0){
\setlength{\unitlength}{0.3\leng} %

\put(20,90){\line(0,1){20}}      %
\put(40,90){\line(0,1){20}}      %
\put(80,90){\line(0,1){20}}      %

\put(20,90){\line(1,0){60}}      %
\put(20,90){\line(1,-1){30}}     %
\put(80,90){\line(-1,-1){30}}    %
\put(50,40){\line(0,1){20}}      %
\put(60,104){\makebox[0pt]{$\cdots$}} }}

\newcommand{\scalecr}[5]{

\put(0,0){
\setunit{0.3\leng} %
\scalecq
\put(20,114){\makebox[0pt]{#1}}   %
\put(40,114){\makebox[0pt]{#2}}   %
\put(80,114){\makebox[0pt]{#3}}   %
\put(50,76){\makebox(0,0){#5}}      
\put(45,38){\makebox(0,0)[t]{#4}}  }} 

\newcommand{\scalecs}[9]{

\put(0,0){
\setunit{0.3\leng} %
\put(20,0){ \scalecq
\put(20,114){\makebox[0pt]{#4}}   %
\put(40,114){\makebox[0pt][r]{#5}}   %
\put(80,114){\makebox[0pt]{#6}}   %
\put(50,76){\makebox(0,0){#9}}      
\put(50,36){\makebox(0,0)[t]{#7}}} \setunit{0.2\leng}
\put(39.8,124){\setleng{0.67\leng} \scalecr{#1}{#2}{#3}{}{#8}}}}

\newcommand{\scalect}[4]{
\setlength{\unitlength}{0.3\leng} %

\put(20,0){\line(0,1){20}} %
\put(20,20){\line(-2,3){20}} %
\put(20,20){\line(2,3){20}}  %
\put(0,50){\line(1,0){40}}  %
\put(0,50){\line(0,1){20}}  %
\put(40,50){\line(0,1){20}} %

\put(15,0){\makebox(0,0)[r]{#3}} %
\put(0,75){\makebox(0,0)[b]{#1}}  %
\put(40,75){\makebox(0,0)[b]{#2}}  %
\put(20,40){\makebox(0,0)[c]{#4}}}

\newcommand{\scalecu}[3]{
\setlength{\unitlength}{0.3\leng} %

\put(20,0){\line(0,1){20}} %
\put(20,20){\line(-2,3){20}} %
\put(20,20){\line(2,3){20}}  %
\put(0,50){\line(1,0){40}}  %
\put(20,50){\line(0,1){20}} %

\put(15,0){\makebox(0,0)[r]{#2}} %
\put(20,75){\makebox(0,0)[b]{#1}}  %
\put(20,40){\makebox(0,0)[c]{#3}}}

\newcommand{\scalecv}[6]{
\put(3,0){\scalect{}{#3}{#4}{#6}} %
\put(-1.5,21){\setleng{0.75\leng} \scalecu{#1}{#2}{#5}}}

\newcommand{\scalecy}[7]{
\put(-9,0){\scalecr{#1}{#2}{#3}{#4}{#6}} %
\put(1.5,-2.5){\setleng{0.75\leng} \scalecu{}{#5}{#7}}}

\usepackage{amsfonts}
\usepackage{amssymb}
\usepackage{amsmath}
\usepackage{fancyheadings}
\usepackage{latexcad}

\begin{document}

\title{Opetopic bicategories: comparison with the classical theory}
\author{Eugenia Cheng\\ \\Department of Pure Mathematics, University
of Cambridge\\E-mail: e.cheng@dpmms.cam.ac.uk}
\date{October 2002}
\maketitle

\begin{abstract}

We continue our previous modifications of the Baez-Dolan theory of
opetopes to modify the Baez-Dolan definition of universality, and
thereby the category of opetopic $n$-categories and lax functors.
For the case $n=2$ we exhibit an equivalence between this category and the category
of bicategories and lax functors.  We examine notions of strictness in the
opetopic theory.

\end{abstract}

\setcounter{tocdepth}{3}
\tableofcontents

\section*{Introduction}
\addcontentsline{toc}{section}{Introduction}

The aim of this paper is to shed some light on the opetopic
definition of weak $n$-category by examining the case $n$=2.  Our
previous work has been on the relationship between different
approaches to the theory of $n$-categories.  In this work we make
a further gesture towards comparison, demonstrating an equivalence
between the opetopic and classical approaches to bicategories.

The opetopic definition we use here is a modified version of the
one given by Baez and Dolan in \cite{bd1}.  Their definition
proceeds in two stages.  First, a language for describing
$k$-cells is set up. Then a concept of universality is introduced,
to deal with composition and coherence.

Our previous work has focused on the construction of $k$-cell
shapes, that is, the theory of opetopes. In \cite{che7} and
\cite{che8} we show that the approach of Baez and Dolan is
equivalent to those of Hermida/Makkai/Power (\cite{hmp1}) and
Leinster (\cite{lei1}), but only with a crucial modification to
their definition.  The key difference is the use of symmetric
multicategories with a {\em category} (rather than a set) of
objects.  We refer the reader to (\cite{che7}) for the full
details.

The definition of $n$-category that we use in this work includes
the above modification, and we refer to the notion thus defined as
`opetopic $n$-category'.

Any proposed definition of $n$-category should
at least be in some way equivalent to the classical definitions as
far as the latter are understood.  We exhibit such equivalence for
the cases $n \leq 2$.

In \cite{che9} we followed through the effects of our previous modifications to
include the recursive definition of opetopic set.  In the present paper we
begin, in Section~\ref{newdef}, by completing this process, to modify the
definition of universality and hence of $n$-category.  The structure of the
definitions is exactly as given in \cite{bd1}; we do not seek to propose a new
approach.  Thus, universality is defined in terms of factorisation properties as in \cite{bd1}.  An opetopic
$n$-category is then defined to be an $I$-opetopic set in which every niche has
a universal occupant, and composites of universals are universal. By this
point, the words of the definition are exactly the same, as all the differences
have been absorbed in the earlier constructions.

Note that in this setting, composites of cells are not necessarily
uniquely defined.  This is one of the key ways in which the theory
differs from the classical theory.

In \cite{bd1} an $n$-functor is defined to be a morphism of the
underlying $I$-opetopic sets, preserving universality.  However,
we consider the more general notion of `lax $n$-functor', in which
universality is not required to be preserved; questions of
strictness are discussed later. So we define the category
    \[\cat{Opic-$n$-Cat}\]
of opetopic $n$-categories and lax functors.  Baez and Dolan do
not construct an $(n+1)$-category of $n$-categories, and we do not
seek to construct one in this work.  Although this leaves the
theory of opetopic $n$-categories still incomplete, we will see
that for $n=2$ a comparison with the classical theory is possible
even without such a construction.

Finally in this section, we restate some useful results from
\cite{bd1}, with the above modifications.

In Section~\ref{prelim} we begin to unravel the definitions for
simple cases in which the interlocking recursion terminates
quickly.  First we examine the cases $n=0$ and $n=1$. We exhibit
equivalences
    \[\cat{Opic-$0$-cat} \simeq \cat{Set}\]
and
    \[\cat{Opic-$1$-cat} \simeq \cat{Cat}.\]

To construct a category \cl{C} from an opetopic 1-category $X$,
the general idea is:
\begin{itemize}
\item the objects of \cl{C} are the 0-cells of $X$
\item the arrows of \cl{C} are the 1-cells of $X$
\item composition is defined by 2-ary universal 2-cells in $X$
\item identities are given by 0-ary universal 2-cells in $X$
\item axioms are seen to hold by considering universal 3-cells in
$X$.
\end{itemize}
This begins to give us a general flavour of how the comparison
proceeds for higher dimensions.

We then discuss some of the properties of $n$-cells in an
$n$-category, which will be useful later.

In Section~\ref{bicatbicat} we examine the case $n=2$ and prove the
main theorem of this work, giving an equivalence
    \[\cat{Opic-$2$-cat} \simeq \cat{Bicat}\]
where \cat{Bicat} is the category of bicategories and lax
functors.

In comparing the opetopic and classical approaches to bicategories
there are two main issues.

\renewcommand{\labelenumi}{\arabic{enumi})}
\begin{enumerate}
\item An opetopic 2-category has $m$-ary 2-cells for all $m \geq
0$, that is, a 2-cell may have a string of $m$ composable 1-cells
as its domain; however a 2-cell in a bicategory has only one
1-cell as its domain.
\item An opetopic 2-category does not have chosen universal
2-cells, that is, 1-cell composition is not uniquely defined;
however, in a bicategory $m$-fold composition is uniquely defined
for $m=0,2$ (identities are considered as 0-fold composites).
\end{enumerate}

The first matter is dealt with in a relatively straightforward
(albeit technically tedious) way, by generating the necessary sets
of $m$-cells from 1-ary 2-cells and 1-cell composites.

The second point involves an issue of {\em choice}.  To construct
a bicategory \cl{B} from an opetopic 2-category $X$ we must make
some choices to give the 0- and 2-fold composites.  The general
idea is:

\begin{itemize}
\item the 0-cells of \cl{B} are the 0-cells of $X$
\item the 1-cells of \cl{B} are the 1-cells of $X$
\item the 2-cells of \cl{B} are the 1-ary 2-cells of $X$
\end{itemize}
\noindent We then make certain choices of 0-ary and 2-ary
universal 2-cells.  Then
\begin{itemize}
\item 1-cell composition in \cl{B} is given by the chosen 2-ary
universal 2-cells
\item 1-cell identities in \cl{B} are given by the chosen nullary
universal 2-cells
\item constraints are induced from composites of the chosen
universals
\item axioms are seen to hold by examining 4-cells.
\end{itemize}

In fact, we define a category $\cat{Opic-$2$-Cat}_b$ of `biased
opetopic 2-categories' whose objects are opetopic 2-categories
equipped with the above choices, but whose morphisms are not
required to preserve those choices.  So the morphisms are simply
morphisms of the underlying 2-categories, and we clearly have an
equivalence
    \[\cat{Opic-$2$-Cat}_b \simeq \cat{Opic-$2$-Cat}.\]
We then exhibit an equivalence
    \[\cat{Opic-$2$-Cat}_b \simeq \cat{Bicat},\]
deferring the more involved calculations to the Appendix.

Finally, in Section~\ref{bicatstr}, we study notions of strictness in the
opetopic theory.  We demonstrate that, while the definition of `lax
$n$-functor' strictifies easily to `weak $n$-functor' and `strict $n$-functor',
the definition of `weak $n$-category' neither laxifies nor strictifies easily.

\subsection*{Terminology}

Although we have modified concepts from \cite{bd1} we have
generally not changed the terminology or notation.

Thus, a \sm\ $Q$ has a category of objects $o(Q)$; an arrow $f$
has source and target $s(f)$ and $t(f)$ respectively.  The slice
$Q^+$ of $Q$ is as defined in \cite{che7} rather than \cite{bd1}.
Baez and Dolan show how to construct a pullback multicategory
given a set over $o(Q)$.  Since with our modification $o(Q)$ is
now a category, we may construct a pullback given a category
\bb{A} and functor $\bb{A} \lra o(Q)$.  We refer the reader to
\cite{che7} for the full definitions.


\subsection*{Acknowledgements}

This work was supported by a PhD grant from EPSRC.  I would like to thank
Martin Hyland and Tom Leinster for their support and guidance.

\label{ncats}

\section{Definitions}
\label{newdef}

In \cite{bd1}, weak $n$-categories are defined as opetopic sets
satisfying certain universality conditions.  In our previous work
(\cite{che7}, \cite{che8}, \cite{che9}) we have examined the
theory of opetopes and opetopic sets.  It remains to discuss the
notion of universality.

We briefly recall here that an opetopic set is a presheaf over the
category of opetopes.  We may think of this as a set of $k$-cells for
each $k \geq 0$, where a $k$-cell is a $k$-opetope with all
constituent $j$-opetopes labelled by $j$-cells in a way that respects
sources and targets.  For the full definition of opetopic sets we
refer the reader to \cite{che9}; for some examples of $k$-cells for
$k \leq 3$ see Section~\ref{exone}.  


\subsection{Universality} \label{bicatuniv}

In the definition of opetopic $n$-category, it is universality
that deals with composition, constraints, axioms and coherence. We
now modify the Baez-Dolan definition of universality in the
context of the results of our earlier work.
Furthermore, with clarity in mind we state the definition in a
terser form than in \cite{bd1}.

\numroman In Section~\ref{bicatncat} we will have the following
definition: {\em An {\em opetopic $n$-category} is an opetopic set
in which
    \begin{enumerate}
    \item Every niche has an $n$-universal occupant.
    \item Every composite of $n$-universals is $n$-universal.
    \end{enumerate}}

\noindent We use the word `composite' in the following sense.  Let
$a_1, \ldots, a_r$ and $c$ be $k$-cells in an opetopic set $X$, with $k \geq
1$. Given a universal $(k+1)$-cell 
	\[u:(a_1, \ldots, a_r) \lra c\] 
we say that $c$ is a composite of $a_1, \ldots, a_r$.  In particular given a
universal cell
	\[u:(a,b) \lra c\] we say that
$c$ is a composite of $a$ and $b$ and also that $u$
and $b$ give a factorisation of $c$ through $a$ (Similarly $u$ and
$a$ give a factorisation of $c$ through $b$).

If $a$ and $b$ are pasted at the target of $b$, say, we may
represent this as
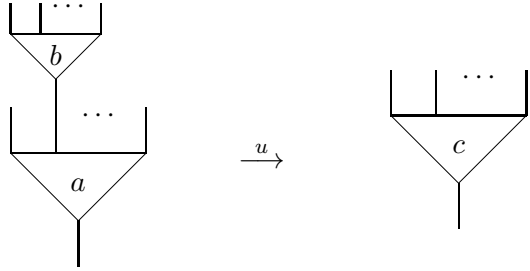
\begin{center} \setleng{1mm}
\begin{picture}(18,37)(8,12)
\scalecs{}{}{}{}{}{}{}{$b$}{$a$}
\end{picture}
\begin{picture}(42,18) \put(15,13){$\map{u}$} \end{picture}
\begin{picture}(18,25)(14,7)
\scalecr{}{}{}{}{$c$}
\end{picture}.
\end{center} Alternatively, regarding $a$, $b$ and $c$ as objects of a
symmetric multicategory at the next dimension up, we may represent
this as

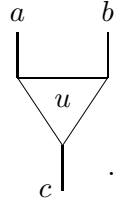
\begin{center} \setleng{1mm}
\begin{picture}(12,25)(0,2)
\scalect{$a$}{$b$}{$c$}{$u$}
\end{picture}.
\end{center}

We now define $n$-universality for $k$-cells and for $k$-cell
factorisations.  The definition is by descending induction on $k$.  Recall that
a niche may be regarded as a potential domain for a cell; so a niche for a
$k$-cell is a $(k-1)$-pasting diagram. We say a cell is `in'
a particular niche if it does indeed have this pasting diagram as its domain.

\begin{definition} A $k$-cell $\alpha$ is {\em $n$-universal} if either
$k>n$ and $\alpha$ is unique in its niche, or $k \leq n$ and (1)
and (2) below are satisfied:

\renewcommand{\labelenumi}{(\arabic{enumi})}
\begin{enumerate}
\item Given any $k$-cell $\gamma$ in the same niche as $\alpha$, there is a
factorisation $u:(\beta,\alpha) \lra \gamma$

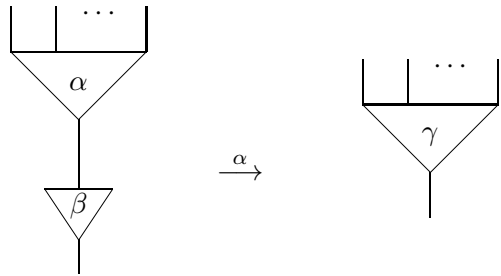
\begin{center}
\setleng{1mm}
\begin{picture}(18,40)(0,-2)
\scalecy{}{}{}{}{}{$\alpha$}{$\beta$}
\end{picture}
\begin{picture}(25,18) \put(5,12){$\map{\alpha}$} \end{picture}
\begin{picture}(18,30)(8,5)
\scalecr{}{}{}{}{$\gamma$}
\end{picture}.
\end{center}

\item Any such factorisation is $n$-universal.

\end{enumerate}
\end{definition}

\begin{definition} A factorisation $u:(b,a) \lra c$ of $k$-cells
is {\em $n$-universal} if $k>n$, or $k\leq n$ and (1) and (2)
below are satisfied:

\renewcommand{\labelenumi}{(\arabic{enumi})}
\begin{enumerate}
\item Given any $k$-cell $b'$ in the same frame as $b$, and any
$(k+1)$-cell \[v:(b',a) \lra c\] with $b'$ and $a$ pasted in the
same configuration as $b$ and $a$ in the source of $u$, there is a
factorisation of $(k+1)$-cells $(u,y)\lra v$

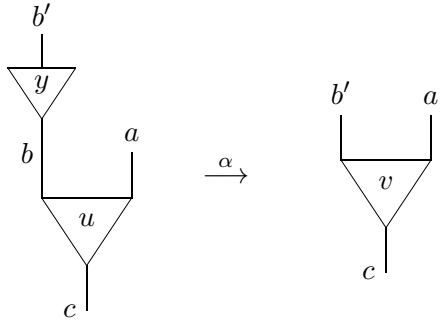
\begin{center}
\setleng{1mm}
\begin{picture}(18,45)(0,0)
\scalecv{$b'$}{$b$}{$a$}{$c$}{$y$}{$u$}
\end{picture}
\begin{picture}(30,18) \put(5,17){$\map{\alpha}$} \end{picture}
\begin{picture}(18,30)(8,-5)
\scalect{$b'$}{$a$}{$c$}{$v$}
\end{picture}
\end{center}

\item Any such factorisation is itself $n$-universal.

\end{enumerate}
\end{definition}
\renewcommand{\labelenumi}{\arabic{enumi})}

If $n$ is clear from the context then we simply say `universal'.

Note that in the terminology of \cite{bd1}, the definition of
`universal factorisation' given above corresponds to a special
case of `balanced punctured niche'.  Furthermore, in each of the
above definitions, each clause (1) and (2) corresponds to the
assertion that a certain punctured niche is balanced.

Although we have still only defined `opetopic $n$-category' in
passing, the following examples concerning particular cases in
opetopic $n$-categories may help to clarify the above definitions.

\begin{examples}
\

{\em \label{bicatcellex}

\begin{enumerate}

\item In an opetopic $n$-category the (unique) universal 1-ary
($n+1$)-cells have the form $x \lra x$, since we have such
universals given by the targets of universal nullary ($n+2$)-cells
    \[(\cdot) \lra (x \rightarrow x).\]

\item In an opetopic $n$-category, a factorisation of $n$-cells is
universal if and only if it is unique.  To see this, consider such
a universal factorisation $u:(b,a) \lra c$.  Now any $(n+1)$-cell
is unique in its niche and hence universal, so any $(n+1)$-cell
$v:(b',a) \lra c$ is a factorisation.  But then, by universality
of the first factorisation, we have a (necessarily universal)
$(n+1)$-cell $y:b' \lra b$ giving $b=b'$ and $u=v$, i.e. the
factorisation is unique.

\item In a 1-category, a 1-cell $x \map{f} y$ is universal if and
only if for any 1-cell $x \map{g} z$ there is a {\em unique}
factorisation

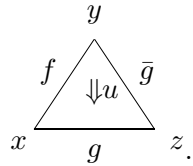
\begin{center}
\begin{picture}(25,22)(-5,-4)
\twotwo{$x$}{$y$}{$z$}{$f$}{$\bar{g}$}{$g$}{$u$}
\end{picture}.
\end{center}

\item In a 2-category, a 1-cell $x \map{f} y$ is universal if and
only if for any 1-cell $x \map{g} z$ there is a factorisation as
above; however, we do not demand that such a factorisation be
unique, but only universal.  That is, given a 2-cell

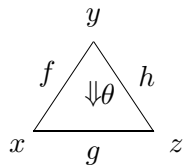
\begin{center}
\begin{picture}(25,22)(-5,-4)
\twotwo{$x$}{$y$}{$z$}{$f$}{$h$}{$g$}{$\theta$}
\end{picture}
\end{center}there is a {\em unique} factorisation

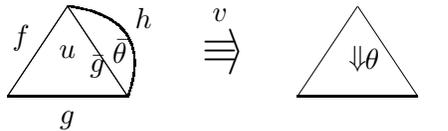
\begin{center}
\begin{picture}(25,22)(-5,-4)
\twotwob{}{}{}{$f$}{}{$g$}{}
\diagbt{}{}{}{}{$h$}{$\bar{g}$}{}{$\bar{\theta}$}{$u$}
\end{picture}
\begin{picture}(11,18) \three{7}{10}{1mm}{$v$} \end{picture}
\begin{picture}(25,22)(-5,-4)
\twotwo{}{}{}{}{}{}{$\theta$}
\end{picture}.
\end{center}

\item In a 3-category, $f$ as above is 3-universal if and only if
any such factorisation $v$ as above is universal (rather than
unique). That is, given any 3-cell

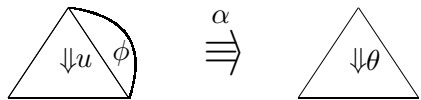
\begin{center}
\begin{picture}(25,22)(-5,-4)
\twotwo{}{}{}{}{}{}{$u$} \diagbt{}{}{}{}{}{}{}{$\phi$}{}
\end{picture}
\begin{picture}(11,18) \three{7}{10}{1mm}{$\alpha$} \end{picture}
\begin{picture}(25,22)(-5,-4)
\twotwo{}{}{}{}{}{}{$\theta$}
\end{picture}
\end{center}

\noindent there is a {\em unique} factorisation
\end{enumerate}}

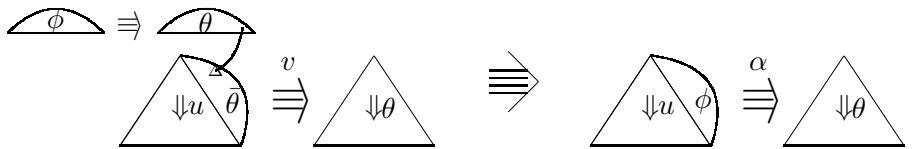
\begin{center}
\begin{picture}(30,28)(-10,-4) \put(-15,15){\diagco{0.8mm}{$\phi$}}
\diagp{1}{16}{4mm} \put(5,15){\diagco{0.8mm}{$\bar{\theta}$}}
\diagr{12}{10}{1mm} \twotwo{}{}{}{}{}{}{$u$}
\diagbt{}{}{}{}{}{}{}{$\bar{\theta}$}{}
\end{picture}
\begin{picture}(6,18) \three{1}{10}{1mm}{$v$} \end{picture}
\begin{picture}(15,22)(3,-4)
\twotwo{}{}{}{}{}{}{$\theta$}
\end{picture}
\begin{picture}(11,18) \diagc{3}{10}{1mm}{} \end{picture}
\begin{picture}(25,22)(-5,-4)
\twotwo{}{}{}{}{}{}{$u$} \diagbt{}{}{}{}{}{}{}{$\phi$}{}
\end{picture}
\begin{picture}(6,18) \three{1}{10}{1mm}{$\alpha$} \end{picture}
\begin{picture}(15,22)(3,-4)
\twotwo{}{}{}{}{}{}{$\theta$}
\end{picture}.
\end{center}

\end{examples}


\renewcommand{\labelenumi}{\roman{enumi})}

\begin{definitions}
\

\begin{itemize}
\item An $n$-{\em
coherent} $Q$-{\em algebra} is a $Q$-opetopic set in which
\begin{enumerate}
\item Every niche has a universal cell in it (or universal
`occupant').
\item Composites of universals are universal.
\end{enumerate}
\item A {\em morphism of $n$-coherent $Q$-algebras} is simply a
morphism of their underlying $Q$-opetopic sets.
\end{itemize}
\end{definitions}

Observe that an $n$-coherent $Q$-algebra is specified uniquely up
to \iso\ by the sets $X(k)$ and functions $f_k$ for $k \leq n+1$,
since for $k \geq n+2$ the sets $X(k)$ and functions $f_k$ are
induced.  A morphism of such is then uniquely determined by the
functions $F_k$ for $k \leq n$.

In \cite{bd1} a morphism of $n$-coherent $Q$-algebras is required
to preserve universality, yielding a stronger notion.  We will
later see that for $n=2$ this gives weak rather than lax functors
of bicategories.  For the time being we consider the lax case
only; we discuss strictness in Section~\ref{bicatstr}.

\subsection{Opetopic $n$-categories} \label{bicatncat}

We are now ready to state the definition of $n$-category.  The
statement here is exactly as in \cite{bd1}; the differences have
all been absorbed into the preliminary definitions.  However, we
note that the exact relationship between our complete modified
definition and the exact Baez-Dolan original remains unclear.

\begin{definitions} \nopagebreak \ \label{ncatdef}
\begin{itemize} \nopagebreak
\item An {\em opetopic $n$-category} is an $n$-coherent $I$-algebra.
\nopagebreak \item A {\em lax $n$-functor} is a morphism of
$n$-coherent $I$-algebras.
\end{itemize} \nopagebreak
We write {\em{\bf Opic-$n$-Cat}} for the category of opetopic
$n$-categories and lax $n$-functors.

\end{definitions}

\renewcommand{\labelenumi}{\roman{enumi})}
\noindent So an opetopic $n$-category is an opetopic set in which
    \begin{enumerate}
    \item Every niche has an $n$-universal occupant.
    \item Every composite of $n$-universals is $n$-universal.
    \end{enumerate}

We now restate, in this modified context, a useful proposition
from \cite{bd1}.  This is a generalisation of the fact that in a
category \cl{C}, for any objects $a,b$ we have a `homset'
$\cl{C}(a,b)$ of morphisms $a \lra b$.  Similarly, in a bicategory
\cl{B}, we have `hom-categories' $\cl{B}(a,b)$ whose objects are
1-cells and morphisms 2-cells; so we also have, for any 2-cells
$\alpha, \beta$, homsets $\cl{B}(\alpha, \beta)$.

Thus in an $n$-category we expect to have `hom-$(n-m)$-categories'
of $m$-cells.  However, since here the domain of an $m$-cell is
not necessarily just a single $(m-1)$-cell, instead of having just
a pair of $(m-1)$-cells as above, we need an $m$-frame to give the
domain and codomain specifying the hom-category.  (Recall that an $m$-frame
consists of an $(m-1)$-pasting diagram together with an $(m-1)$-cell that might
be the domain and codomain of an $m$-cell.)

\begin{proposition} \label{bicatpropb1} Let $X$ be an $n$-coherent $Q$-algebra.
Then for $m \leq n$ any $m$-frame determines an opetopic
$(n-m)$-category.
\end{proposition}

The idea is first to restrict $X$ to cells of dimension $m$ and
above; this is clearly still $(n-m)$-coherent.  We can then
restrict to only those cells in the given frame $\alpha$ by
`pulling back' along the morphism
    \[1 \map{\alpha} S(m).\]
So we follow Baez-Dolan and use the following construction of
`pullback opetopic set'. Let $Q$ and $Q'$ be tidy \sms\ with
object-categories \bb{C} and $\bb{C}'$ respectively, with $\bb{C}
\simeq S$ and $\bb{C'} \simeq S'$ discrete.  Let $X$ be a
$Q$-opetopic set. Suppose we have a morphism $S' \longrightarrow
S$.  Then we may construct a pullback opetopic set $X'$ by
induction as follows. Let $X'(0)$ be given by the pullback

\setlength{\unitlength}{0.2em}
\begin{center}
\begin{picture}(45,40)(5,5)  %

\put(10,10){\makebox(0,0){$S'(0)$}}  
\put(10,35){\makebox(0,0){$X'(0)$}}  
\put(45,35){\makebox(0,0){$X(0)$}}  
\put(45,10){\makebox(0,0){$S(0)$}}  

\put(17,35){\vector(1,0){21}}  
\put(17,10){\vector(1,0){21}}  
\put(10,30){\vector(0,-1){15}} 
\put(45,30){\vector(0,-1){15}} 

\put(52,9){\makebox(0,0){.}} 

\end{picture}
\end{center} Now we have equivalences
    \[\begin{array}{c} o({Q_{{X(0)}}}^+) \simarr S(1),\\
    o({Q'_{{X(0)}}}^+) \simarr S'(1) \end{array}\]
where $S(1)$ and $S'(1)$ are discrete.  So the morphism
    \[X'(0) \lra X(0)\]
induces a morphism
    \[S'(1) \lra S(1)\]
and we may form a pullback opetopic set of $X_1$ along this
morphism; we set this to be $X'_1$, the underlying
${Q'_{{X'(0)}}}^+$-opetopic set of $X'$.

\begin{proposition} \label{bicatpropa1} {\em(See \cite{bd1}, Proposition
45)} If $X$ is $n$-coherent then $X'$ is $n$-coherent.
\end{proposition}

\begin{prf} It is easy to check that a cell in $X'$ is universal
if and only if the corresponding cell in $X$ is universal, and
that a factorisation in $X'$ is universal if and only if the
corresponding factorisation in $X$ is universal.  \end{prf}

\begin{prfof}{Proposition \ref{bicatpropb1}} Let $\alpha$ be an
$m$-frame in $X$ with $m \leq n$, so $\alpha \in S(m)$.  Now $X$
determines an $(n-m)$-coherent $Q(m)$-algebra, and we have a
morphism
    \[o(I)=1 \map{\alpha} S(m)\]
so we may form a pullback $I$-opetopic set along this morphism.

By Proposition~\ref{bicatpropa1} this is $(n-m)$-coherent, i.e. it
is an opetopic $(n-m)$-category.  \end{prfof}

\begin{examples}  \label{bicat1cat} \end{examples}
\begin{enumerate}\renewcommand{\labelenumi}{\arabic{enumi})}
\item In an $n$-category $X$, every 1-frame determines an
$(n-1)$-category.

A 1-frame in $X$ is given by

\begin{center}
\begin{picture}(22,8)(2,0) \onecell{0}{0}{1mm}{$a$}{$b$}{$?$}
\end{picture}\hspace{5mm}.
\end{center} We denote the induced
$(n-1)$-category by Hom$(a,b)$ or $X(a,b)$; its cells are of the
form shown below.

\noindent 0-cells \sunit
\begin{picture}(35,10)(-25,0)
\put(5,0){\vector(1,0){24}}  %
\put(4,0){\makebox(0,0)[r]{$a$}} %
\put(31,0){\makebox(0,0)[l]{$b$}} %
\put(17,1){\makebox(0,0)[b]{$f$}}
\end{picture}

\noindent 1-cells
\begin{picture}(20,20)(-15,0) \onetwo{15}{0}{1mm}{$a$}{$b$}{$f$}{$g$}{$\alpha$}
\end{picture}

\noindent 2-cells ($k$-ary) \sunit
\begin{picture}(40,25)(-15,0)
\setlength{\unitlength}{6mm}
\put(0,0){\line(1,0){4}}  %
\qbezier(0,0)(2,1.2)(4,0)   %
\qbezier(0,0)(2,3)(4,0) %
\qbezier(0,0)(2,4.2)(4,0) %
\put(2.1,0.2){\makebox(0,0)[c]{$\alpha_k$}}
\put(2,1.2){\makebox(0,0)[c]{$\vdots$}}
\put(2.1,1.8){\makebox(0,0)[c]{$\alpha_1$}}

\diagj{7}{0}{6mm}         %
\diagp{5.5}{0.5}{9mm}    %
\put(5.5,1.2){\makebox(0,0)[b]{$\theta$}}  %
\put(9,0.5){\makebox(0,0)[c]{$\alpha$}}

\end{picture}

\vspace{1cm}\hspace*{5mm}\vdots

\item Given a 2-frame

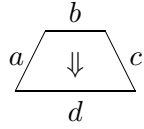
\sunit
\begin{picture}(20,20)(-25,-5)
\setlength{\unitlength}{0.4mm}

\put(0,0){\line(1,2){10}}  %
\put(10,20){\line(1,0){20}}  %
\put(30,20){\line(1,-2){10}}  %
\put(0,0){\line(1,0){40}}

\put(3,11){\makebox(0,0)[r]{$a$}}  %
\put(20,23){\makebox(0,0)[b]{$b$}}  %
\put(20,-3){\makebox(0,0)[t]{$d$}}  %
\put(38,11){\makebox(0,0)[l]{$c$}}  %
\put(20,9){\makebox(0,0)[c]{$\Downarrow$}}  %

\end{picture}

\noindent say, we have an $(n-2)$-category whose cells are of the
form shown below.

\noindent 0-cells
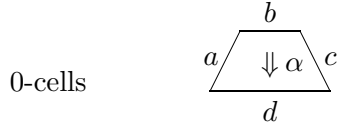
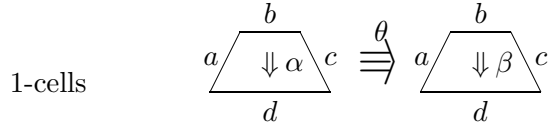
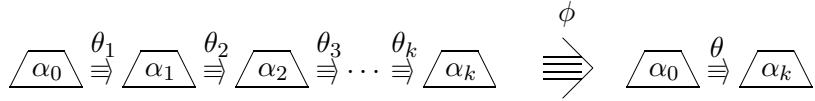
\begin{picture}(20,20)(-15,0)
\setlength{\unitlength}{0.4mm}

\put(0,0){\line(1,2){10}}  %
\put(10,20){\line(1,0){20}}  %
\put(30,20){\line(1,-2){10}}  %
\put(0,0){\line(1,0){40}}

\put(3,11){\makebox(0,0)[r]{$a$}}  %
\put(20,23){\makebox(0,0)[b]{$b$}}  %
\put(20,-3){\makebox(0,0)[t]{$d$}}  %
\put(38,11){\makebox(0,0)[l]{$c$}}  %
\put(20,9){\makebox(0,0)[c]{$\Downarrow$}}  %
\put(25,9){\makebox(0,0)[l]{$\alpha$}}

\end{picture}

\noindent 1-cells \sunit
\begin{picture}(40,20)(-15,0)
\setlength{\unitlength}{0.4mm}

\put(0,0){\begin{picture}(20,20)

\put(0,0){\line(1,2){10}}  %
\put(10,20){\line(1,0){20}}  %
\put(30,20){\line(1,-2){10}}  %
\put(0,0){\line(1,0){40}}

\put(3,11){\makebox(0,0)[r]{$a$}}  %
\put(20,23){\makebox(0,0)[b]{$b$}}  %
\put(20,-3){\makebox(0,0)[t]{$d$}}  %
\put(38,11){\makebox(0,0)[l]{$c$}}  %
\put(20,9){\makebox(0,0)[c]{$\Downarrow$}}  %
\put(25,9){\makebox(0,0)[l]{$\alpha$}} \end{picture}}

\put(70,0){\begin{picture}(20,20)

\put(0,0){\line(1,2){10}}  %
\put(10,20){\line(1,0){20}}  %
\put(30,20){\line(1,-2){10}}  %
\put(0,0){\line(1,0){40}}

\put(3,11){\makebox(0,0)[r]{$a$}}  %
\put(20,23){\makebox(0,0)[b]{$b$}}  %
\put(20,-3){\makebox(0,0)[t]{$d$}}  %
\put(38,11){\makebox(0,0)[l]{$c$}}  %
\put(20,9){\makebox(0,0)[c]{$\Downarrow$}}  %
\put(25,9){\makebox(0,0)[l]{$\beta$}}
\end{picture}
}

\diagp{55}{10}{6mm}  %
\put(57,17){\makebox(0,0)[b]{$\theta$}}

\end{picture}

\vspace{1cm} \noindent 2-cells ($k$-ary)\nopagebreak

\begin{picture}(100,10) %
\setlength{\leng}{1mm} \setlength{\unitlength}{\leng}
\diagbl{0}{0}{0.6\leng}  %
\put(5,2){\makebox(0,0)[c]{$\alpha_0$}} %
\diagbl{15}{0}{0.6\leng}  %
\put(20,2){\makebox(0,0)[c]{$\alpha_1$}} %
\diagbl{30}{0}{0.6\leng}  %
\put(35,2){\makebox(0,0)[c]{$\alpha_2$}} %
\diagp{12}{2}{3.6\leng}  %
\put(12.5,4){\makebox(0,0)[b]{$\theta_1$}} %
\diagp{27}{2}{3.6\leng}  %
\put(27.5,4){\makebox(0,0)[b]{$\theta_2$}} %
\diagp{42}{2}{3.6\leng}  %
\put(42.5,4){\makebox(0,0)[b]{$\theta_3$}} %
\put(47,2){\makebox(0,0)[c]{$\ldots$}}  %
\diagp{52}{2}{3.6\leng}  %
\put(52.5,4){\makebox(0,0)[b]{$\theta_k$}} %
\diagbl{55}{0}{0.6\leng}  %
\put(60,2){\makebox(0,0)[c]{$\alpha_k$}} %
\diagc{70}{0}{\leng}  %
\put(74,8){\makebox(0,0)[b]{$\phi$}} %
\diagbl{82}{0}{0.6\leng} %
\put(87,2){\makebox(0,0)[c]{$\alpha_0$}} %
\diagp{94}{2}{3.6\leng}  %
\put(94,4){\makebox(0,0)[b]{$\theta$}} %
\diagbl{97}{0}{0.6\leng} %
\put(102,2){\makebox(0,0)[c]{$\alpha_k$}} %

\end{picture}

 \vspace{1cm}\hspace*{5mm}\vdots

\item Given an $(n-1)$-frame we have a 1-category whose
objects are $(n-1)$-cells and arrows are 1-ary $n$-cells.

\end{enumerate}

\section{Preliminary examples} 
\label{prelim}

Any proposed definition of $n$-category should at least be in some
way equivalent to the classical definitions as far as the latter are
understood.  In \cite{bd1} Baez and Dolan examine the case $n=1$ but
do not explain how their definition is equivalent to the classical
definition of bicategories in the case $n=2$.  This is perhaps
because, without the modifications described in this our earlier
work, such an equivalence does not arise.  The main results of this
work gives an equivalence between the (modified) opetopic and the
classical approaches to bicategories.  We begin in this section with
some examples to help clarify and motivate the later arguments; our
general aim is to shed some light on the inescapable loops in the
definition of universality, as well as to compare the resulting
structures with the classical ones. 

Note that for $n \leq 1$ the difference between our definition and
the original Baez-Dolan definition is not yet apparent.  The
result for $n=1$ is described in \cite{bd1} (Example 42); we
include it here (with more detail) for completeness.

\subsection{Opetopic 0-categories} \label{exzero}

An opetopic 0-category $X$ is determined, up to \iso, by the set
$X(0)$.  For, given any 0-cell $a \in A$, the following nullary
2-niche \length{0.2em}

\begin{center}
\begin{picture}(30,15) \nulltwo{0}{0}{1mm}{$a$}{$$}{$$} \end{picture}
\end{center}

\noindent must have a unique occupant, and so the unique occupant
of the following 1-niche

\begin{center}
\begin{picture}(30,8) \onecell{0}{0}{1mm}{$a$}{$?$}{$?$} \end{picture}
\end{center}

\noindent must have $a$ as its target, and we can call the 1-cell
$1_a$, giving
    \[X(1) \cong \{a \lra a : a \in A\}.\]

\begin{proposition} There is an equivalence
    \[\mbox{{\bf Opic-$0$-Cat}} \simarr \mbox{\upshape{\bfseries Set}}\]
surjective in the direction shown.
\end{proposition}

\begin{prf} We construct such a functor, $\zeta$. Let $X$ be an opetopic
0-category. We put \[\zeta(X) = X(0).\] A morphism $f:X \lra Y$ of
opetopic 0-categories is uniquely specified by the function
$f_0:X(0) \lra Y(0)$ so we put \[\zeta(f) = f_0.\]

\noindent Conversely, given a set $A$, we have an opetopic
0-category $X$ such that $\zeta(X)=A$; $X$ is defined by
    \begin{eqnarray*} X(0) &=& A\\
    X(1) &=& \{a \map{1_a} a : a \in A\}.\end{eqnarray*}
So $\zeta$ is surjective, and it is clearly full and faithful,
giving an equivalence as required. \end{prf}

\subsection{Opetopic 1-categories} \label{exone}

We first clarify our notation.  We draw

\begin{itemize}
\item 1-cells as arrows
\begin{picture}(20,10) \onecell{0}{0}{1mm}{}{}{} \end{picture}
\item 2-cells as
\begin{picture}(20,30) \diagbp{18}{0}{1mm} \end{picture}
\item 3-cells as
\begin{picture}(100,30) \diagbo{18}{0}{1mm} \end{picture} $\ldots$
\end{itemize}

\noindent These represent \iso\ classes of objects in the
appropriate \sm. We give below some typical examples of openings,
niches, frames and cells.

\begin{picture}(100,15)
\put(0,0){\makebox(0,0)[l]{1-opening}}  %
\onecell{50}{0}{1mm}{}{}{}
\end{picture}

\begin{picture}(100,15)
\put(0,0){\makebox(0,0)[l]{1-niche}}  %
\onecell{50}{0}{1mm}{$a$}{?}{?}
\end{picture}

\begin{picture}(100,15)
\put(0,0){\makebox(0,0)[l]{1-frame}}  %
\onecell{50}{0}{1mm}{$a$}{$b$}{?}
\end{picture}

\begin{picture}(100,15)
\put(0,0){\makebox(0,0)[l]{1-cell}}  %
\onecell{50}{0}{1mm}{$a$}{$b$}{$f$}
\end{picture}

\vspace{2cm} \hspace*{4.5cm} 3-ary \hspace{4.2cm} nullary

\begin{picture}(100,30)
\put(0,0){\makebox(0,0)[l]{2-opening}}  %
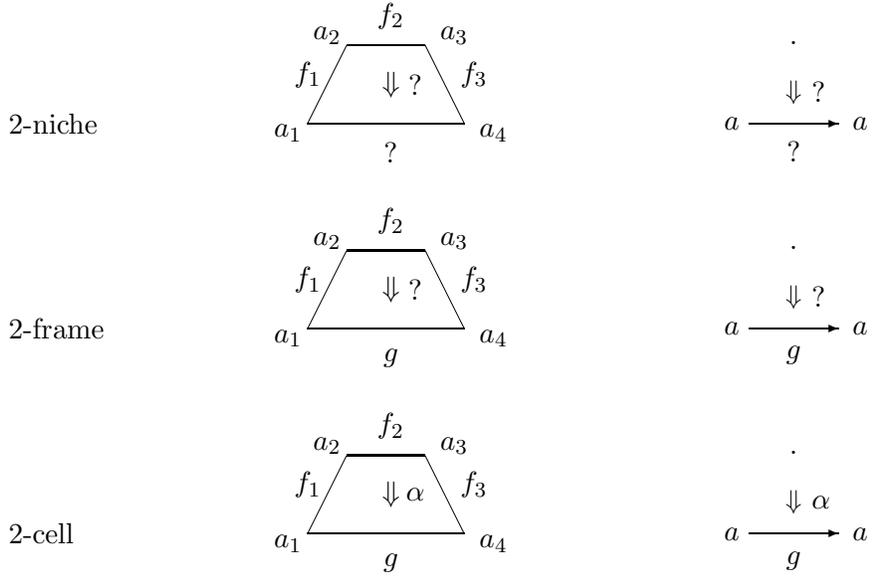
\put(50,0){
\begin{picture}(0,0)
\threetwo{$a_1$}{$a_2$}{$a_3$}{$a_4$}{$$}{$$}{$$}{$$}{$$}
\end{picture}}
\nulltwo{120}{0}{1mm}{$a$}{}{}
\end{picture}

\begin{picture}(100,35)
\put(0,0){\makebox(0,0)[l]{2-niche}}  %
\put(50,0){
\begin{picture}(0,0)
\threetwo{$a_1$}{$a_2$}{$a_3$}{$a_4$}{$f_1$}{$f_2$}{$f_3$}{?}{?}
\end{picture}}
\nulltwo{120}{0}{1mm}{$a$}{?}{?}
\end{picture}

\begin{picture}(100,35)
\put(0,0){\makebox(0,0)[l]{2-frame}}  %
\put(50,0){
\begin{picture}(0,0)
\threetwo{$a_1$}{$a_2$}{$a_3$}{$a_4$}{$f_1$}{$f_2$}{$f_3$}{$g$}{?}
\end{picture}}
\nulltwo{120}{0}{1mm}{$a$}{$g$}{?}
\end{picture}

\begin{picture}(100,45)(0,-10)
\put(0,0){\makebox(0,0)[l]{2-cell}}  %
\put(50,0){
\begin{picture}(0,0)
\threetwo{$a_1$}{$a_2$}{$a_3$}{$a_4$}{$f_1$}{$f_2$}{$f_3$}{$g$}{$\alpha$}
\end{picture}}
\nulltwo{120}{0}{1mm}{$a$}{$g$}{$\alpha$}
\end{picture}

\noindent Where confusion is unlikely, we may omit some
lower-dimensional labels once the higher-dimensional ones are in
place, as in the following examples.

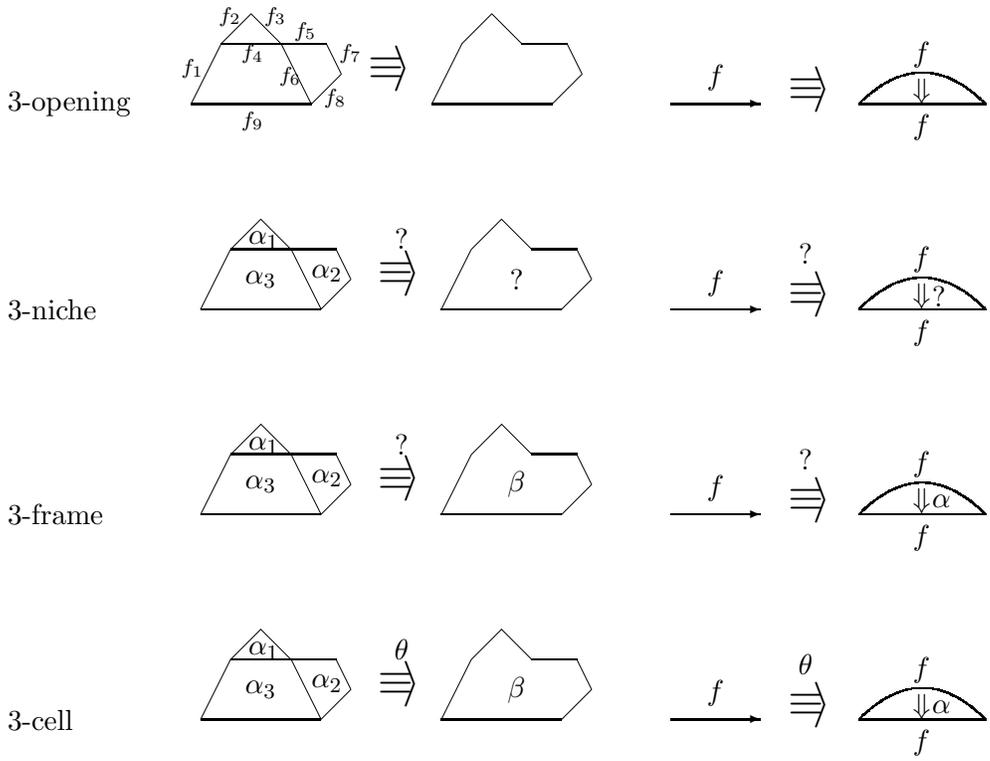
\begin{picture}(100,30)(10,0)
\put(0,0){\makebox(0,0)[l]{3-opening}}  %
\put(30,0){\begin{picture}(50,50) \fontsize{8pt}{15} \sunit
\diagbo{0}{0}{1mm} %

\put(0,4){$f_1$} %
\put(5,11){$f_2$} %
\put(11,11){$f_3$} %
\put(8,6){$f_4$} %
\put(15,9){$f_5$} %
\put(13,3){$f_6$} %
\put(21,6){$f_7$} %
\put(19,0){$f_8$} %
\put(8,-3){$f_9$} %

\end{picture}}
\put(105,0){ \begin{picture}(0,0) \nullthree{$f$}{}{}
\end{picture}}
\end{picture}

\begin{picture}(100,35)(10,0)
\put(0,0){\makebox(0,0)[l]{3-niche}}  %
\put(30,0){
\begin{picture}(0,0)
\threethree{$\alpha_1$}{$\alpha_2$}{$\alpha_3$}{?}{?}
\end{picture}}
\put(105,0){ \begin{picture}(0,0) \nullthree{$f$}{?}{?}
\end{picture}}
\end{picture}

\begin{picture}(100,35)(10,0)
\put(0,0){\makebox(0,0)[l]{3-frame}}  %
\put(30,0){
\begin{picture}(0,0)
\threethree{$\alpha_1$}{$\alpha_2$}{$\alpha_3$}{$\beta$}{?}
\end{picture}}
\put(105,0){ \begin{picture}(0,0) \nullthree{$f$}{$\alpha$}{?}
\end{picture}}
\end{picture}

\begin{picture}(100,35)(10,0)
\put(0,0){\makebox(0,0)[l]{3-cell}}  %
\put(30,0){
\begin{picture}(0,0)
\threethree{$\alpha_1$}{$\alpha_2$}{$\alpha_3$}{$\beta$}{$\theta$}
\end{picture}}
\put(105,0){ \begin{picture}(0,0)
\nullthree{$f$}{$\alpha$}{$\theta$} \end{picture}}
\end{picture}


\subsubsection*{}

We begin by constructing a functor
    \[\zeta : \mbox{{\bf Opic-$1$-Cat}} \lra \cat{Cat};\]

\noindent we will eventually show that this functor is an
equivalence.

\begin{itemize}
\item On objects
\end{itemize}

Given an opetopic 1-category $X$ we define a category
$\cl{C}=\cl{C}_X$ as follows. First set ob \cl{C} $=X(0)$. Then,
given objects $a,b \in X(0)$, let $\cl{C}(a,b)$ be the preimage of
$a \map{?} b$ under $f_1$. (Recall that we have a 0-category
Hom$(a,b)$, that is, a set.)

Composition and identities in \cl{C} are defined according to the
2-cells in $X$ as follows.  For composition consider 1-cells $a
\map{f} b$, $b \map{g} c$.  We have the following 2-niche

\begin{center}
\begin{picture}(25,22)(-5,-4)
\twotwo{}{}{}{$f$}{$g$}{?}{?}
\end{picture}
\end{center}

\noindent which has a unique occupant; we write it as

\begin{center}
\begin{picture}(25,22)(-5,-4)
\twotwo{}{}{}{$f$}{$g$}{$gf$}{$u$}
\end{picture}\ \ .
\end{center}

\noindent  For identities we have already observed (Examples
\ref{bicatcellex}) that in an opetopic $n$-category the universal
1-ary $(n+1)$-cells are of the form $a \lra a$.  Explicitly, for
$n=1$ we have for any $a \in X(0)$ a nullary 2-niche

\begin{center}
\begin{picture}(25,18)(0,-5)
\nulltwo{0}{0}{1mm}{$a$}{$?$}{$?$}
\end{picture}
\end{center}

\noindent which must have a unique occupant. So we write it as

\begin{center}
\begin{picture}(25,18)(0,-5)
\nulltwo{0}{0}{1mm}{$a$}{$1_a$}{$u$}
\end{picture}
\end{center}

\noindent  and check that this does indeed act as the identity
with respect to the composition defined above.  We seek the unique
occupant of the niche

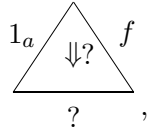
\begin{center}
\begin{picture}(27,22)(-5,-4)
\twotwo{}{}{}{$1_a$}{$f$}{?}{?}
\end{picture},
\end{center}

\noindent that is

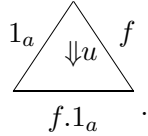
\begin{center}
\begin{picture}(27,22)(-5,-4)
\twotwo{}{}{}{$1_a$}{$f$}{$f.1_a$}{$u$}
\end{picture}.
\end{center}

\noindent Certainly we have the following 3-niche

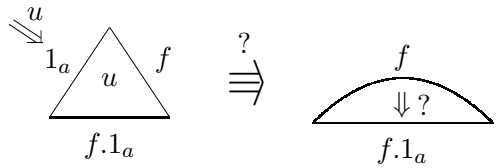
\begin{center}
\begin{picture}(30,30)(-9,-4)
\diagbw{$$}{$$}{$$}{$1_a$}{$f$}{$f.1_a$}{$u$}{$u$}
\end{picture}
\begin{picture}(16,18) \three{11}{10}{1mm}{?} \end{picture}
\begin{picture}(35,15)(5,-3)
\onetwo{10}{0}{1mm}{$$}{$$}{$f$}{$f.1_a$}{$?$}
\end{picture}
\end{center}

\noindent with a unique occupant.  So by
Example~\ref{bicatcellex}(1), we have $f.1_a = f$ as required.
Similarly $1_a.f = f$.

It remains to check that associativity holds.  Given 1-cells
    \[a \map{f} b \map{g} c \map{h} d\]
we have the following universal 3-cells

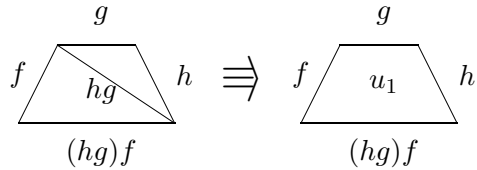
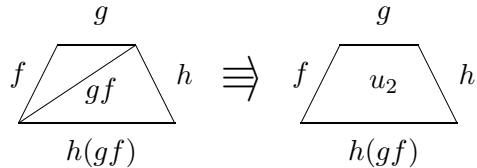
\sunit

\begin{center}
\begin{picture}(25,20)(-2,-4)
\assleft{$f$}{$g$}{$h$}{$hg$}{$(hg)f$}
\end{picture}
\begin{picture}(10,18) \diagp{5}{10}{6mm} \end{picture}
\begin{picture}(25,20)(-2,-4)
\threetwob{$$}{$$}{$$}{$$}{$f$}{$g$}{$h$}{$(hg)f$}{$u_1$}
\end{picture}
\end{center}

\begin{center}
\begin{picture}(25,30)(-2,-4)
\assright{$f$}{$g$}{$h$}{$gf$}{$h(gf)$}
\end{picture}
\begin{picture}(10,18) \diagp{5}{10}{6mm} \end{picture}
\begin{picture}(25,20)(-2,-4)
\threetwob{$$}{$$}{$$}{$$}{$f$}{$g$}{$h$}{$h(gf)$}{$u_2$}
\end{picture}
\end{center}

\noindent  But $u_1$ and $u_2$ are occupants of the same 2-niche;
by uniqueness they must be the same, giving
    \[(hg)f = h(gf)\]
as required.  So we have defined a category, and we set
    \[\zeta(X) = \cl{C}_X.\]
Observe that we find composites and identities by considering
universal 2-cells, and we check axioms by considering universal
3-cells.

\begin{itemize}
\item On morphisms
\end{itemize}

Given a morphism of opetopic 1-categories $F:X \lra Y$ we seek to
define a functor $F:\cl{C}_X \lra \cl{C}_Y$.  We define the action
of $F$ on objects and arrows by the functions
    \begin{eqnarray*} F_0 &:& X(0) \lra Y(0)\\
    \mbox{and \ \ } F_1 &:& X(1) \lra Y(1).\end{eqnarray*}
We check functoriality.  By definition of morphisms of opetopic
1-categories, the following diagram commutes

\setlength{\unitlength}{0.2em}
\begin{center}
\begin{picture}(45,40)(10,5)      %

\put(10,10){\makebox(0,0){$Y(1)$}}  
\put(10,35){\makebox(0,0){$X(1)$}}  
\put(45,35){\makebox(0,0)[l]{$o({I_{X(0)}}^+) = o(Q(1))$}}  %
\put(45,10){\makebox(0,0)[l]{$o({I_{Y(0)}}^+) = o(R(1))$}}  %

\put(18,35){\vector(1,0){25}}  
\put(18,10){\vector(1,0){25}}  
\put(10,30){\vector(0,-1){15}} 
\put(53,30){\vector(0,-1){15}} 

\put(8,23){\makebox(0,0)[r]{$F_1$}} 

\end{picture}
\end{center} giving
    \begin{eqnarray*} F(\mbox{dom }f) &=& \mbox{dom (Ff)}\\
    \mbox{and \ \ } F(\mbox{cod }f)  &=& \mbox{cod (Ff)}.\end{eqnarray*}
Now the function
    \[F_2:X(2) \lra Y(2)\]
makes the following diagram commute

\setlength{\unitlength}{0.2em}
\begin{center}
\begin{picture}(45,40)(5,5)      %

\put(10,10){\makebox(0,0){$Y(2)$}}  
\put(10,35){\makebox(0,0){$X(2)$}}  
\put(45,35){\makebox(0,0)[l]{$o({Q(1)_{X(1)}}^+)$}}  %
\put(45,10){\makebox(0,0)[l]{$o({R(1)_{Y(1)}}^+)$}}  %

\put(18,35){\vector(1,0){25}}  
\put(18,10){\vector(1,0){25}}  
\put(10,30){\vector(0,-1){15}} 
\put(53,30){\vector(0,-1){15}} 

\put(8,23){\makebox(0,0)[r]{$F_2$}} 
\put(85,8){\makebox(0,0){}}

\end{picture}
\end{center} so under the action of $F_2$ the following (universal) 2-cell in $X$
\sunit
\begin{center}
\begin{picture}(25,22)(-5,-4)
\twotwo{}{}{}{$f$}{$g$}{$gf$}{$u$}
\end{picture}
\end{center}

\noindent gives the following 2-cell in $Y$

\sunit
\begin{center}
\begin{picture}(25,22)(-5,-4)
\twotwob{}{}{}{$Ff$}{$Fg$}{$F(gf)$}{$Fu$}
\end{picture}
\end{center}

\noindent and so we have $F(g \circ f) = Fg \circ Ff$ by
uniqueness of 2-niche occupants.  Similarly consider the following
nullary 2-cell in $X$

\begin{center}
\begin{picture}(25,18)(0,-5)
\nulltwo{0}{0}{1mm}{$a$}{$1_a$}{$u$}
\end{picture}.
\end{center}


\noindent Under the action of $F_2$ we have the following 2-cell
in $Y$

\begin{center}
\begin{picture}(25,18)(0,-5)
\nulltwo{0}{0}{1mm}{$Fa$}{$F(1_a)$}{$Fu$}
\end{picture}
\end{center}

\noindent and so we have $F(1_a)= 1_{Fa}$ by uniqueness of 2-niche
occupants.

So $F$ is a functor as required. Observe that in the above
construction we do not need to stipulate that universality be
preserved.

\bigskip

Finally, before showing that $\zeta$ is an equivalence, we
characterise universal 1-cells as invertibles.

\begin{proposition} \label{bicatyon} A 1-cell $f$ in $X$ is universal
if and only if it is invertible as an arrow of $\cl{C}_X$.
\end{proposition}

\noindent {\bfseries Proof 1 (bare hands).} Let $a \map{f} b$ be a
universal 1-cell in $X$.  We certainly have a 1-cell

\begin{center}
\begin{picture}(22,8)(2,0) \onecell{0}{0}{1mm}{$a$}{$a$}{$1_a$}
\end{picture}.
\end{center} So by clause (1) of the definition of universal 1-cell we have a
factorisation, that is a 1-cell
\begin{center}
\begin{picture}(22,8)(2,0) \onecell{0}{0}{1mm}{$b$}{$a$}{$g$}
\end{picture}
\end{center}and a universal 2-cell

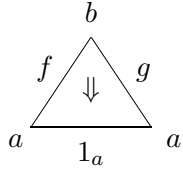
\begin{center}
\begin{picture}(25,22)(-3.5,-4)
\twotwo{$a$}{$b$}{$a$}{$f$}{$g$}{$1_a$}{}
\end{picture}
\end{center}

\noindent so we have $gf = 1_a$.

Now consider the 1-cell
\begin{center}
\begin{picture}(22,8)(2,0) \onecell{0}{0}{1mm}{$a$}{$b$}{$f$}
\end{picture}.
\end{center}Similarly, we have a universal 2-cell

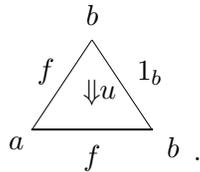
\begin{center}
\begin{picture}(25,22)(-3.5,-4)
\twotwo{$a$}{$b$}{$b$}{$f$}{$1_b$}{$f$}{$u$}
\end{picture}.
\end{center}

\noindent Now by clause (1) of the definition of universal 2-cell,
if we have a 2-cell

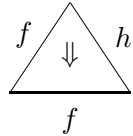
\begin{center}
\begin{picture}(25,22)(-3.5,-4)
\twotwo{$$}{$$}{$$}{$f$}{$h$}{$f$}{$$}
\end{picture}
\end{center}

\noindent then we have a factorisation, so we certainly have  a
2-cell

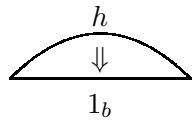
\begin{center}
\begin{picture}(35,15)(5,-3)
\onetwo{10}{0}{1mm}{$$}{$$}{$h$}{$1_b$}{$$}
\end{picture}
\end{center}

\noindent By uniqueness of 2-niche occupants, this gives
    \[hf = f \Rightarrow h = 1_b.\]
Now consider the following 3-cell

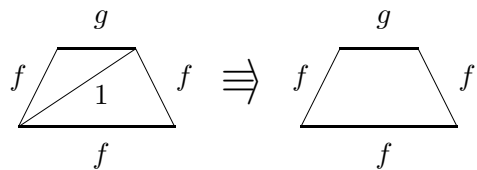
\sunit
\begin{center}
\begin{picture}(25,30)(-2,-4)
\assright{$f$}{$g$}{$f$}{$1$}{$f$}
\end{picture}
\begin{picture}(10,18) \diagp{5}{10}{6mm} \end{picture}
\begin{picture}(25,20)(-2,-4)
\threetwob{$$}{$$}{$$}{$$}{$f$}{$g$}{$f$}{$f$}{$$}
\end{picture}
\end{center}

\noindent giving $f(gf) = f$.  But by associativity we have
    \[f(gf) = (fg)f = f\]
so we have $fg = 1_b$. So if $f$ is universal in $X$ then $f$ is
invertible in $\cl{C}_X$.

Conversely, suppose $f$ is invertible in $\cl{C}_X$, so we have in
$X$ 2-cells

\begin{center}
\begin{picture}(25,27)(-5,-4)
\twotwo{$a$}{$b$}{$a$}{$f$}{$g$}{$1$}{}
\end{picture}
\ \ \ \ ,\
\begin{picture}(25,22)(-5,-4)
\twotwo{$b$}{$a$}{$b$}{$g$}{$f$}{$1$}{}
\end{picture}.
\end{center}

\noindent We now show that $f$ is universal:

\begin{enumerate}
\item Given any 0-cell $b' \in X(0)$ and 1-cell $a \map{h} b'$ we have the following
3-cell

\begin{center}
\begin{picture}(25,22)(-2,-4)
\assright{$f$}{$g$}{$h$}{$1$}{$h$}
\end{picture}
\begin{picture}(10,18) \diagp{5}{10}{6mm} \end{picture}
\begin{picture}(25,20)(-2,-4)
\threetwob{$$}{$$}{$$}{$$}{$f$}{$g$}{$h$}{$h$}{$$}
\end{picture}
\end{center}

\noindent so by associativity the following universal 2-cell

\begin{center}
\begin{picture}(25,22)(-3.5,-4)
\twotwo{}{}{}{$f$}{$hg$}{$h$}{}
\end{picture}
\end{center}

\noindent giving a factorisation for $h$ as required.

\item We show that any such factorisation is universal.  Let

\begin{center}
\begin{picture}(25,22)(-5,-4)
\twotwo{}{}{}{$f$}{$s$}{$h$}{}
\end{picture}
\end{center}

\noindent be such a factorisation.  Then given any other 2-cell

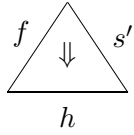
\begin{center}
\begin{picture}(25,22)(-5,-4)
\twotwo{}{}{}{$f$}{$s'$}{$h$}{}
\end{picture}
\end{center}

\noindent we need to exhibit a factorisation

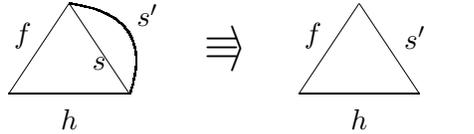
\begin{center}
\begin{picture}(25,22)(-5,-4)
\diagbt{$$}{$$}{$$}{$f$}{$s'$}{$s$}{$h$}{$$}{$$}
\end{picture}
\begin{picture}(11,18) \three{7}{10}{1mm}{} \end{picture}
\begin{picture}(28,22)(-5,-4)
\twotwob{$$}{$$}{$$}{$f$}{$s'$}{$h$}{$$}
\end{picture}.
\end{center}

\noindent Now \[h=sf \Rightarrow hg=sfg=s\] so we have $s'=hg=s$
and 3-cell

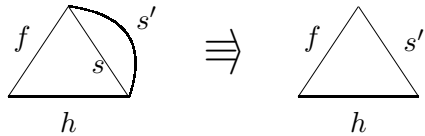
\begin{center}
\begin{picture}(25,22)(-5,-4)
\diagbt{$$}{$$}{$$}{$f$}{$s'$}{$s$}{$h$}{$$}{$$}
\end{picture}
\begin{picture}(11,18) \three{7}{10}{1mm}{} \end{picture}
\begin{picture}(25,22)(-5,-4)
\twotwob{$$}{$$}{$$}{$f$}{$s'$}{$h$}{}
\end{picture}
\end{center} as required.  Any such factorisation is then trivially
universal.
\end{enumerate}

\noindent So if $f$ is invertible then $f$ is universal, and the
proposition is proved. \qeed

\bigskip

Although the above calculations may help in understanding the
definitions, the proposition may be proved more quickly using the
Yoneda Lemma as follows.

\bigskip

\renewcommand{\labelenumi}{\arabic{enumi})}
\noindent {\bfseries Proof 2 (Yoneda).  } $f$ is universal in $X$
if and only if
\begin{enumerate}
\item Given any arrow $b \map{g} c$ there is an arrow $b
\map{\bar{g}} c$ such that $\bar{g}f=g$ and
\item $h_1 f = h_2 f \Rightarrow h_1 = h_2$
\end{enumerate}

\noindent i.e. for all $c \in \mbox{ob } \cl{C}$ the function
    \[\begin{array}{rccc}
    f^\ast : & \cl{C}(b,c) & \lra & \cl{C}(a,c)
    \\
    & h & \longmapsto & h \circ f \\
    \end{array}\]
is an \iso.  But this is true if and only if $f$ is \iso\ since
the Yoneda embedding is full and faithful.  \qeed

\bigskip

In a later work (\cite{che11}) we propose a characterisation of
universality that generalises the above Yoneda result.

\begin{proposition}\label{bicatpropc1} The functor $\zeta$ exhibits an
equivalence of categories
    \[\cat{Opic-$1$-Cat} \simarr \cat{Cat}\]
surjective in the direction shown.
\end{proposition}

\begin{prf} We have defined a functor
    \[\zeta:\cat{Opic-$1$-Cat} \lra \cat{Cat}\]
above, and it is clearly full and faithful; we show that it is
surjective.

Given any (small) category \cl{C}, we may construct an opetopic
1-category $X$ with $X(0)=\mbox{ob }\cl{C}$ and $X(1)=\mbox{arr
}\cl{C}$. We see immediately that every 1-niche has a universal
occupant $a \map{1_a} a$.  The set $X(2)$ is defined as follows.

Every nullary 2-niche

\begin{center}
\begin{picture}(25,18)(0,-5)
\nulltwo{0}{0}{1mm}{$a$}{$?$}{$?$}
\end{picture}
\end{center}

\noindent has a unique occupant

\begin{center}
\begin{picture}(25,18)(0,-5)
\nulltwo{0}{0}{1mm}{$a$}{$1_a$}{$$}
\end{picture}
\end{center}

\noindent and every $m$-ary 2-niche

\begin{center}
\begin{picture}(25,25)(-2,-5)
\mtwo{$f_1$}{$f_2$}{$f_m$}{$?$}{$?$}
\end{picture}
\end{center}

\noindent has a unique occupant

\begin{center}
\begin{picture}(25,25)(-2,-5)
\mtwo{$f_1$}{$f_2$}{$f_m$}{$f_m \circ f_{m-1} \circ \ldots \circ
f_1$}{}
\end{picture}.
\end{center} Furthermore, since a 1-cell is universal if and only if
it is invertible as an arrow of \cl{C}, composites of universals
are universal.

So $X$ is 1-coherent, and clearly $\zeta(X) = \cl{C}$.
\end{prf}


\subsection{$n$-cells in an $n$-category} \label{bicatsecd1}


The definition of universality works from the top down: universal
cells are understood via cells in the dimension above, and the
starting point is that all cells in dimensions higher than $(n+1)$
are trivial.  So in effect, $n$-cells result from the `first' step
of the induction; we now make some general observations about
$n$-cells, which will be useful later.

Recall (Example \ref{bicat1cat}(3)) that every $(n-1)$-frame
determines an opetopic 1-category.  So we have an opetopic
1-category of $(n-1)$-cells and 1-ary $n$-cells, or, by
Proposition~\ref{bicatpropc1}, a category.

Let $X$ be an opetopic $n$-category.  First recall that composites
of $n$-cells in $X$ are uniquely determined, since occupants of
$(n+1)$-niches are unique.  Also, composition of $n$-cells is
strictly associative and a morphism of opetopic $n$-categories
must be strictly functorial on $n$-cell composites.  (In fact, we
have a \sm\ of $(n-1)$-cells and $n$-cells.)

Now consider an $n$-niche $\alpha$ in $X$.  Then, given any
universal occupant $u$, every occupant $f$ of $\alpha$ factors
uniquely as
    \[f=g \circ u\]
where g is a 1-ary $n$-cell.  So, for any such universal, we may
express the set of occupants of $\alpha$ as
    \[g \circ u \mbox{ \ such that \ } g \in X(n)_1 \mbox{ and } s(g) = t(u)\]
where $X(n)_1$ is the set of 1-ary $n$-cells. Given any other
universal occupant $u'$, we then have
    \[u' = x \circ u\]
for some (unique) universal $x$.  So we have
    \[\{g' \circ u' \} = \{g \circ u\}\]
since $g' \circ u' = g' \circ (x \circ u) = (g' \circ x) \circ u$.

More generally, given any non-empty set $U$ of universal occupants
of $\alpha$, the set of occupants of $\alpha$ may be expressed as
    \[ \begin{array} {r} \{g \circ u : u \in U, g \in X(n)_1,
    s(g) = t(u)\} \\ \  \end{array}
    \left/ \begin{array}{l} \\ \sim \end{array} \right.  .\]
Here $\sim$ is the equivalence relation generated by
\begin{enumerate}
\item $g \circ u \sim g' \circ u' \iff g = g' \circ x_{uu'}$
\item $1 \circ u \sim u$
\end{enumerate}
where for any $u,u' \in U$, $x_{uu'}$ is the unique universal such
that \[u' = x_{uu'} \circ u.\]

\section{Bicategories} 
\label{bicatbicat}

\label{bicatex}

We are now ready to turn our attention to the case $n=2$.  We show
how to construct a classical bicategory from an opetopic 2-category,
leading to the main theorem, which shows how the opetopic and
classical theories of bicategories are equivalent.

An important difference between this construction and that for the
case $n=1$ is that an element of choice now arises.  The
universality condition stipulates that every niche should have a
universal occupant, but does not {\em specify} such universals.
This approach differs from the approach of Leinster
({\cite{lei5}}), for example, in which contractibility is defined
as a property but specific contractions are then given.


This approach also differs from the classical approach to
bicategories, in which binary and nullary composites of 1-cells
are specified, even though $m$-fold composites are not, for $m>2$.
(Note that 1-cell identities are considered as `nullary
composites'.)  Leinster refers to this theory as being `biased'
towards binary composites; in {\cite{lei1}}, he introduces the
notion of {\em unbiased bicategory}. The theory of bicategories is
made `unbiased' by specifying $m$-fold composites for all $m$.
This theory turns out to be equivalent to the classical one
(\cite{lei5}). Leinster also comments that, provided at least one
choice has been made for each of $k=0$ and some $k\geq 2$, an
equivalent theory of bicategories may be formed.

Another way of eliminating bias from a bicategory might be to
choose {\em no specified composites}.  We will later see that this
is how the opetopic approach may be interpreted. Once we have
shown that this theory is equivalent to the classical one, it is
easy to see which choices give rise to a theory of bicategories,
and it follows immediately that all such theories are equivalent. 
This issue turns out to be related to the question of {\em
strictness}, and we discuss these notions in Section~\ref{bicatstr}.

\subsection{The main theorem: equivalence with the classical theory}

We show that the opetopic and classical theories of bicategories are
equivalent, in the following sense.

\begin{theorem}\label{bicatthmd1}  Write \cat{Bicat} for the category
of bicategories and morphisms (lax functors).  Then
    \[\cat{Opic-$2$-Cat} \simeq \cat{Bicat}.\]
\end{theorem}

Given an opetopic 2-category X, we seek to construct a bicategory
\cl{B} (using the definition given in \cite{lei2}).  To do this we
need to make some choices of universal 2-cells.  The general idea
is

\begin{itemize}
\item the 0-cells of \cl{B} are the 0-cells of $X$
\item the 1-cells of \cl{B} are the 1-cells of $X$
\item the 2-cells of \cl{B} are the 1-ary 2-cells of $X$.
\end{itemize}
\noindent We then choose a universal occupant for each 0-ary and
2-ary 2-niche in $X$.  Then
\begin{itemize}
\item 1-cell composition in \cl{B} is given by the chosen 2-ary
universal 2-cells in $X$
\item 1-cell identities in \cl{B} are given by the chosen nullary
universal 2-cells in $X$
\item constraints are induced from composites of the chosen
universals
\item axioms are seen to hold by examining 4-cells.
\end{itemize}

\noindent In fact, we define a category of `biased opetopic
2-categories' in which these choices have already been made.

\begin{definitions} \
\begin{itemize}
\item A {\em biased opetopic 2-category} is an opetopic 2-category
together with a chosen universal occupant for every nullary and
2-ary 2-niche.
\item A {\em morphism of biased opetopic 2-categories} is simply a
morphism of the underlying 2-categories.
\end{itemize}
We write $\cat{Opic-$2$-Cat}_b$ for the category of biased
opetopic 2-categories and morphisms.
\end{definitions}

Note that the choice of universal 2-cells is free, that is, the
chosen cells are not required to satisfy any axioms.  Furthermore,
no preservation condition is imposed on the morphisms in this
category.

\begin{proposition} \label{bicatpropf1} There is an equivalence
    \[\cat{Opic-$2$-Cat}_b \simarr \cat{Opic-$2$-Cat}\]
surjective in the direction shown.
\end{proposition}

\begin{prf} Clear from the definitions. \end{prf}


\noindent So in fact, we prove the following proposition:

\begin{proposition} \label{bicatprope1} There is an equivalence
    \[\cat{Opic-$2$-Cat}_b \simarr \cat{Bicat}\]
surjective in the direction shown.
\end{proposition}

\noindent Finally we will make some comments about the choices
made in forming a biased opetopic 2-category.

For the longer calculations in this subsection, and for an
explanation of the `shorthand' used in manipulating 2-cells, we
refer the reader to Appendix C.

\begin{prfof}{Proposition \ref{bicatprope1}} We construct a functor
    \[\zeta:\cat{Opic-$2$-Cat}_b \lra \cat{Bicat}\]
and show that it is surjective, full and faithful.

\begin{itemize} \item We define the action of $\zeta$ on objects.
\end{itemize}

Let $X$ be a biased opetopic 2-category.  So in addition to the
usual data, we have

\renewcommand{\labelenumi}{\roman{enumi})}
\begin{enumerate}
\item for each object $A \in X(0)$ a chosen universal 2-cell

\begin{center}
\begin{picture}(25,18)(0,-5)
\nulltwo{0}{0}{1mm}{$A$}{$$}{$\iota_A$}
\end{picture}
\end{center}

\item for each pair $f,g$ of composable 1-cells, a chosen
universal 2-cell

\begin{center}
\begin{picture}(25,22)(-5,-4)
\twotwob{}{}{}{$f$}{$g$}{}{$\ c_{gf}$}
\end{picture}.
\end{center}
\end{enumerate}

\noindent We may indicate these chosen 2-cells by $\sim$ as in

\begin{center}
\begin{picture}(25,18)(0,-5)
\nulltwo{0}{0}{1mm}{$A$}{}{$\sim$}
\end{picture}
\ \ \  ,
\begin{picture}(29,22)(-5,-4)
\twotwob{}{}{}{$f$}{$g$}{}{$\sim$}
\end{picture}.
\end{center} We now define a bicategory $\cl{B}=\cl{B}_X$ as follows.  First
set
    \[\mbox{ob}(\cl{B}) = X(0).\]
Recall (Proposition \ref{bicatpropb1}) that given objects $A,B \in
X(0)$, we have an opetopic 1-category Hom$(A,B)$. Let
$\cl{B}(A,B)$ be the category corresponding to Hom$(A,B)$
according to Proposition \ref{bicatpropc1}. So we have 1-cells
given by 1-cells of $X$
    \[a \map{f} b\]
and 2-cells given by 1-ary 2-cells of $X$

\begin{center}
\begin{picture}(35,15)(5,-3)
\onetwo{10}{0}{1mm}{$$}{$$}{$f$}{$g$}{$\alpha$}
\end{picture}.
\end{center}

\noindent 2-cell composites are given by the (unique) 3-cell
occupants, for example

\begin{center}
\begin{picture}(25,15)(-5,-4)
\diagco{1mm}{{\small $\alpha$}} %
\diagby{1mm}{{\small $\beta$}}%
\end{picture}
\begin{picture}(11,18) \three{7}{5}{1mm}{} \end{picture}
\begin{picture}(25,15)(-5,-4)
\diagco{1mm}{{\small $\beta \circ \alpha$}}
\end{picture}
\end{center}

\noindent and 2-cell identities by nullary 3-cells

\begin{center}
\begin{picture}(50,15)(3,-5)
\nullthree{$f$}{$1_f$}{$$}
\end{picture}.
\end{center}

Now for any objects $A,B,C \in \mbox{ob } \cl{B}$ we need a
functor
    \[\begin{array}{rccc}
    c_{ABC} : & \cl{B}(B,C) \times \cl{B}(A,B) & \lra & \cl{B}(A,C)
    \\
    & (g,f) & \longmapsto & g \circ f = gf \\
    & (\beta,\alpha) & \longmapsto & \beta \ast \alpha. \\
\end{array}\]
We define $g \circ f$ to be the target 1-cell of the chosen
universal $c_{gf}$, so we have

\begin{center}
\begin{picture}(25,22)(-5,-4)
\twotwob{}{}{}{$f$}{$g$}{$g\circ f$}{$\ c_{gf}$}
\end{picture}.
\end{center}

\noindent Note that for each composable pair $f,g$, we have
specified a 2-cell $c_{gf}$; this is crucially stronger than
merely specifying a 1-cell $g \circ f$.

We now show how horizontal 2-cell composition is induced. Consider
2-cells

\begin{center}
\begin{picture}(35,15)(5,-3)
\onetwo{10}{0}{1mm}{$$}{$$}{$f_1$}{$f_2$}{$\alpha$}
\end{picture}
\ \ ,
\begin{picture}(35,15)(5,-3)
\onetwo{10}{0}{1mm}{$$}{$$}{$g_1$}{$g_2$}{$\beta$}
\end{picture};
\end{center}

\noindent we seek a 2-cell

\begin{center}
\begin{picture}(35,15)(5,-3)
\onetwob{10}{0}{1mm}{$$}{$$}{$g_1.f_1$}{$g_2.f_2$}{$\beta \ast
\alpha$}
\end{picture}.
\end{center}

\noindent We have a 3-cell

\begin{center}

\begin{picture}(30,22)(-5,-4)
\diagbt{$$}{$$}{$$}{$$}{$g_1$}{$g_2\ $}{$g_2f_2$}{$\beta$}{$\sim$}
\diagbv{$$}{$$}{$$}{$f_1$}{$\ f_2$}{$$}{$$}{$\alpha$}{}
\end{picture}
\begin{picture}(7,18) \three{3}{10}{1mm}{} \end{picture}
\begin{picture}(25,22)(-5,-4)
\twotwob{}{}{}{$f_1$}{$g_1$}{$g_1.f_1$}{$\phi$}
\end{picture}
\end{center}
\noindent unique in its niche, and a universal 2-cell

\begin{center}
\begin{picture}(25,22)(-5,-4)
\twotwob{$$}{$$}{$$}{$f_1$}{$g_1$}{$g_1.f_1$}{$\sim$}
\end{picture}
\end{center}

\noindent inducing, by definition of universality, a 2-cell

\begin{center}
\begin{picture}(35,15)(5,-3)
\onetwo{10}{0}{1mm}{$$}{$$}{$g_1.f_1$}{$g_2.f_2$}{$\theta$}
\end{picture}
\end{center}

\noindent unique such that there is a 3-cell

\begin{center}
\begin{picture}(25,22)(-5,-4)
\diagbr{$$}{$$}{$$}{$$}{$$}{$$}{$$}{$\sim$}{$\theta$}
\end{picture}
\begin{picture}(11,18) \three{7}{10}{1mm}{} \end{picture}
\begin{picture}(29,22)(-5,-4)
\twotwob{$$}{$$}{$$}{$$}{$$}{$$}{$\phi$}
\end{picture}.
\end{center}

\noindent Put $\beta \ast \alpha = \theta$.  We check
functoriality, that is
\begin{enumerate}
\item $1_g \ast 1_f = 1_{gf}$
\item $(\beta_2 \circ \beta_1) \ast (\alpha_2 \circ \alpha_1) =
(\beta_2 \ast \alpha_2) \circ (\beta_1 \ast \alpha_1)$ (middle 4
interchange)
\end{enumerate}
(see Appendix, Lemma~\ref{bicatlemma1}).

Next we need, for each object $A$, a 1-cell $A \map{I_A} A$.  We
define this to be the target of the chosen universal $\iota_A$, so
we have

\begin{center}
\begin{picture}(25,18)(0,-5)
\nulltwo{0}{0}{1mm}{$A$}{$I_A$}{$\iota_A$}
\end{picture}.
\end{center}

\noindent Note that, as before, we have specified a universal
2-cell, not just the 1-cell $I_A$.

We now seek natural isomorphisms $a,r,l$.  Each of these is
induced uniquely from the chosen universals $\iota$ and $c$. For
$a$, consider 1-cells
    \[A \map{f} B \map{g} C \map{h} D.\]
We seek a natural isomorphism
    \[a_{hgf} : (hg)f \simarr h(gf).\]
We have

\begin{center}
\begin{picture}(25,23)(-2,-4)
\assleft{$f$}{$g$}{$h$}{$hg$}{$(hg)f$}
\end{picture}
\begin{picture}(10,18) \diagp{5}{10}{6mm} \end{picture}
\begin{picture}(25,20)(-2,-4)
\threetwob{$$}{$$}{$$}{$$}{$f$}{$g$}{$h$}{$(hg)f$}{$\theta$}
\end{picture}
\end{center}and
\begin{center}
\begin{picture}(25,23)(-2,-4)
\assright{$f$}{$g$}{$h$}{$gf$}{$h(gf)$}
\end{picture}
\begin{picture}(10,18) \diagp{5}{10}{6mm} \end{picture}
\begin{picture}(25,20)(-2,-4)
\threetwob{$$}{$$}{$$}{$$}{$f$}{$g$}{$h$}{$h(gf)$}{$\phi$}
\end{picture}.
\end{center}

\noindent $\theta$ and $\phi$ are composites of universals, so
universal.  Universality of $\theta$ induces a unique 2-cell
$\alpha$ such that

\begin{center}
\begin{picture}(25,20)(-2,-4)
\threetwob{}{}{}{}{}{}{}{}{$\theta$} \diagcd{$\alpha$}
\end{picture}
\begin{picture}(12,18) \put(5,7){=} \end{picture}
\begin{picture}(25,20)(-2,-4)
\threetwob{$$}{$$}{$$}{$$}{$$}{$$}{$$}{$$}{$\phi$}
\end{picture}
\end{center}
\noindent so

\begin{center}
\begin{picture}(25,20)(-2,-4)
\assleft{$f$}{$g$}{$h$}{$hg$}{$$} %
\diagcd{$\alpha$}
\end{picture}
\begin{picture}(12,18) \put(5,7){=} \end{picture}
\begin{picture}(25,20)(-2,-4)
\assright{$f$}{$g$}{$h$}{$gf$}{$h(gf)$}
\end{picture}.
\end{center}

\noindent Put $a_{hgf} = \alpha$.  We see from universality of
$\phi$ that $a_{hgf}$ is an \iso; we check that it satisfies
naturality (see Appendix, Lemma~\ref{bicatlemma2}).


Next we seek a natural transformation $r$, so we need for any
1-cell $A \map{f} B$ a 2-cell

\begin{center}
\begin{picture}(35,15)(5,-3)
\onetwo{10}{0}{1mm}{$$}{$$}{$f.I_A$}{$f$}{$r$}
\end{picture}.
\end{center}

\noindent Now we have a 3-cell

\begin{center}
\begin{picture}(30,22)(-9,-4)
\diagbw{$$}{$$}{$$}{$I_A$}{$f$}{$f.I_A$}{$\sim$}{$\sim$}
\end{picture}
\begin{picture}(11,18) \three{7}{10}{1mm}{} \end{picture}
\begin{picture}(35,15)(5,-3)
\onetwo{10}{0}{1mm}{$$}{$$}{$f$}{$f.I_A$}{$\alpha$}
\end{picture}
\end{center}

\noindent and the target 2-cell $\alpha$ is universal since it is
the composite of universals. (Note that this is not the same
$\alpha$ as above.) So $\alpha$ induces

\begin{center}
\begin{picture}(25,10)(-5,-2)
\diagco{1mm}{$\alpha$} %
\diagby{1mm}{$r$}%
\end{picture}
\begin{picture}(11,10) \three{7}{3}{1mm}{} \end{picture}
\begin{picture}(25,10)(-5,-2)
\diagco{1mm}{$1_f$}
\end{picture}
\end{center}

\noindent so

\sunit
\begin{center}
\begin{picture}(30,22)(-9,-4)
\diagbw{$$}{$$}{$$}{$I_A$}{$f$}{$$}{$\sim$}{$$}
\diagbr{$$}{$$}{$$}{$$}{$$}{$$}{$$}{$\sim$}{$r_f$}
\end{picture}
\begin{picture}(12,18) \put(5,7){=} \end{picture}
\begin{picture}(35,15)(5,-3)
\onetwo{10}{0}{1mm}{$$}{$$}{$f$}{$f$}{$1_f$}
\end{picture}.
\end{center}

\noindent Since $\alpha$ is universal it is an \iso\ with $r_f$ as
its inverse; so $r_f$ is also an \iso.  We also check naturality
(see Appendix, Lemma~\ref{bicatlemma3}).  The construction of and
result for $l$ follow similarly.

Finally we check the axioms for a bicategory (see Appendix,
Lemma~\ref{bicatlemma4}). So we have defined a bicategory
$\cl{B}_X$ and we put $\zeta(X)=\cl{B}_X$.

\begin{itemize} \item We define the action of $\zeta$ on
morphisms. \end{itemize}

Let $F:X \lra X'$ be a morphism of opetopic 2-categories, so for
each $k$ we have

\begin{center}
\setlength{\unitlength}{0.2em}
\begin{picture}(48,38)(5,7)     %

\put(10,10){\makebox(0,0){$X'(k)$}}  
\put(10,35){\makebox(0,0){$X(k)$}}  
\put(45,35){\makebox(0,0){$S(k)$}}  
\put(45,10){\makebox(0,0){$S'(k)$}}  

\put(17,35){\vector(1,0){21}}  
\put(17,10){\vector(1,0){21}}  
\put(10,30){\vector(0,-1){15}} 
\put(45,30){\vector(0,-1){15}} 

\put(8,23){\makebox(0,0)[r]{$F_k$}} 
\put(47,23){\makebox(0,0)[l]{}} 
\put(27,37){\makebox(0,0)[b]{$f_k$}} 
\put(27,8){\makebox(0,0)[t]{$f'_k$}} 

\put(55,9){\makebox(0,0){.}} 

\end{picture}
\end{center}

\noindent We construct from $F$ a lax functor
    \[(F,\phi) : \cl{B}_X \lra \cl{B}_{X'}.\]
The action of $F$ on objects is given by the function
    \[F_0: X(0) \lra X'(0);\]
we also need, for any objects $A,B \in \mbox{ob } \cl{B}_X$ a
functor
    \[F_{AB}:\cl{B}_X(A,B) \lra \cl{B}_{X'}(FA,FB).  \]
Now for any $A,B \in \mbox{ob } \cl{B}_X$ we have an opetopic
1-category Hom$(A,B)$, and restricting $F$ to this gives a
morphism of opetopic 1-categories
    \[\mbox{Hom}(A,B) \lra \mbox{Hom}(FA,FB)\]
so by Proposition~\ref{bicatpropc1} we have a functor $F_{AB}$ as
required.

Next we seek a natural transformation $\phi_{ABC}$, so for any
1-cells
    \[A \map{f} B \map{g} C\]
we need a 2-cell
    \[\phi_{gf}: Fg \circ Ff \lra F(g \circ f).\]
We have in $X$ a chosen universal 2-cell

\sunit
\begin{center}
\begin{picture}(25,22)(-5,-4)
\twotwob{}{}{}{$f$}{$g$}{$gf$}{$c$}
\end{picture}
\end{center}
\noindent so under the action of $F$ we have in $X'$ a 2-cell

\begin{center}
\begin{picture}(25,22)(-5,-4)
\twotwob{}{}{}{$Ff$}{$Fg$}{$F(gf)$}{$Fc$}
\end{picture}.
\end{center}
\noindent But in $X'$ we have a chosen universal 2-cell

\begin{center}
\begin{picture}(25,22)(-5,-4)
\twotwob{}{}{}{$Ff$}{$Fg$}{$Fg.Ff$}{$\sim$}
\end{picture}
\end{center}

\noindent which, by definition of universality, induces a 2-cell

\begin{center}
\begin{picture}(35,15)(5,-3)
\onetwob{10}{0}{1mm}{$$}{$$}{$Fg.Ff$}{$F(g.f)$}{$\phi_{gf}$}
\end{picture}
\end{center}

\noindent unique such that

\begin{center}
\begin{picture}(25,22)(-5,-4)
\diagbr{$$}{$$}{$$}{$Ff$}{$Fg$}{$$}{$$}{$$}{$\phi_{gf}$}
\end{picture}
\begin{picture}(12,18) \put(5,7){=} \end{picture}
\begin{picture}(25,22)(-5,-4)
\twotwob{$$}{$$}{$$}{$Ff$}{$Fg$}{$F(g.f)$}{$Fc$}
\end{picture};
\end{center}

\noindent we check that this satisfies naturality (see Appendix,
Lemma~\ref{bicatlemma5}).

We now seek a natural transformation $\phi_A$ for each object $A$,
so we seek a 2-cell

\begin{center}
\begin{picture}(35,15)(5,-3)
\onetwo{10}{0}{1mm}{$$}{$$}{$I_{FA}$}{$FI_A$}{$\phi_A$}
\end{picture}.
\end{center}

\noindent We have in $X$ a chosen universal 2-cell

\begin{center}
\begin{picture}(25,18)(0,-5)
\nulltwo{0}{0}{1mm}{$A$}{$I_A$}{$\iota_A$}
\end{picture}
\end{center}

\noindent so applying $F$ gives a 2-cell in $X'$

\begin{center}
\begin{picture}(25,18)(0,-5)
\nulltwo{0}{0}{1mm}{$FA$}{$FI_A$}{$F\iota_A$}
\end{picture}.
\end{center}

\noindent Now the chosen universal in $X'$

\begin{center}
\begin{picture}(25,18)(0,-5)
\nulltwo{0}{0}{1mm}{$FA$}{$I_{FA}$}{$\iota_{FA}$}
\end{picture}
\end{center}

\noindent induces, by universality, a 2-cell

\begin{center}
\begin{picture}(35,15)(5,-3)
\onetwob{10}{0}{1mm}{$$}{$$}{$I_{FA}$}{$FI_A$}{$\phi_A$}
\end{picture}
\end{center}

\noindent unique such that

\begin{center}
\begin{picture}(25,30)(0,-5)
\diagbq{$$}{$I_{FA}$}{$FI_A$}{$$}{$\phi_A$}
\end{picture}
\begin{picture}(12,18) \put(6,7){=} \end{picture}
\begin{picture}(25,18)(0,-5)
\nulltwo{0}{0}{1mm}{$$}{$FI_A$}{$F\iota_A$}
\end{picture}
\end{center}

\noindent and there is no non-trivial naturality to check.

Finally we check that the axioms for a lax functor hold (see
Appendix, Lemma~\ref{bicatlemma6}).  So $(F,\phi)$ is indeed a lax
functor, and we set $\zeta(F) = (F, \phi)$.

\bigskip

It is clear that the above construction of $\zeta$ is functorial,
so we have defined a functor
    \[\zeta:\cat{Opic-$2$-Cat}_b \lra \cat{Bicat};\]
it remains to show that $\zeta$ is surjective, full and faithful.

\begin{itemize} \item We show that $\zeta$ is surjective.
\end{itemize}

Given a bicategory \cl{B}, we construct an opetopic 2-category $X$
such that $\zeta(X)=\cl{B}$.  The idea is

\renewcommand{\labelenumi}{\roman{enumi})}
\begin{enumerate}
\item The 0-cells of $X$ are the 0-cells of \cl{B}.
\item The 1-cells of $X$ are the 1-cells of \cl{B}.
\item The 1-ary 2-cells of $X$ are the 2-cells of \cl{B}.
\item For $m \neq 1$, certain $m$-ary universals are fixed
according to $m$-fold composites in \cl{B}; the remaining cells
are then generated to ensure that these do indeed satisfy
universality.
\item The 3-cells of $X$ are determined from 2-cell composition in
\cl{B}.
\end{enumerate}

Put $X(0) = \mbox{ob}(\cl{B})$ and set $X(1)$ to be the set of
1-cells of \cl{B}; the function $f_1:X(1) \lra S(1)$ is defined so
that the preimage of the frame $A \map{?} B$ is the set of objects
of the category $\cl{B}(A,B)$.

We now construct X(2) bearing in mind the comments in
Section~\ref{bicatsecd1}.  Write $X(2)_m \subset X(2)$ for the set
of $m$-ary 2-cells.  First we define the set $X(2)_1$ of 1-ary
2-cells to be the set of 2-cells of \cl{B}.

For 0-ary 2-cells, we first define for each $A \in X(0)$ a 2-cell

\begin{center}
\begin{picture}(25,18)(0,-5)
\nulltwo{0}{0}{1mm}{$A$}{$I_A$}{$\iota_A$}
\end{picture}.
\end{center}

\noindent We then define the set of occupants of the same niche to
be
    \[\{\alpha \circ \iota_A : \alpha \in X(2)_1, s(\alpha) =
    I_A\}\]
that is, cells of the form

\begin{center}
\begin{picture}(25,18)(0,-5)
\nulltwo{0}{0}{1mm}{$A$}{$f$}{$\alpha\circ \iota_A$}
\end{picture}
\begin{picture}(12,18) \put(5,7){=} \end{picture}
\begin{picture}(25,30)(0,-5)
\diagbq{$$}{$$}{$$}{$\iota_A$}{$\alpha$}
\end{picture}
\end{center}

\noindent where we put $1 \circ \iota = \iota$.

Similarly for $X(2)_2$ we first define for each composable pair of
1-cells $f,g$ a 2-cell

\begin{center}
\begin{picture}(25,15)(-5,-1)
\twotwob{}{}{}{$f$}{$g$}{$gf$}{$\ c_{gf}$}
\end{picture}
\end{center}
\noindent where $g \circ f$ is the composite in \cl{B}.  We then
define the set of occupants of this niche to be
    \[\{\alpha \circ c_{gf} : \alpha \in X(2)_1, s(\alpha) = g
    \circ f\}\]
that is, cells of the form

\begin{center}
\begin{picture}(25,15)(-5,-1)
\diagbr{$$}{$$}{$$}{$f$}{$g$}{$$}{$$}{$c_{gf}$}{$\alpha$}
\end{picture}
\end{center}

\noindent where we put $1 \circ c = c$.

For $X(2)_m$, $m>2$, consider a 2-niche of the form

\begin{center}
\begin{picture}(25,25)(-2,-5)
\mtwo{$f_1$}{$f_2$}{$f_m$}{$f$}{$?$}
\end{picture}.
\end{center}

\noindent  We have no preferred $m$-fold composite in \cl{B};
instead, for each composite $\gamma(f_1, \ldots, f_m)$ we define a
2-cell $u_\gamma$ which is to be universal:
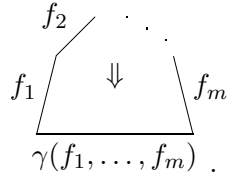
\begin{center}
\begin{picture}(25,25)(-2,-5)
\mtwo{$f_1$}{$f_2$}{$f_m$}{$\gamma(f_1, \ldots, f_m)$}{}
\end{picture}.
\end{center}

\noindent  Now, suppose we have composites $\gamma(f_1, \ldots,
f_m)$ and $\gamma'(f_1, \ldots, f_m)$.  Then we have a unique
invertible
    \[a_{\gamma \gamma'} : \gamma(f_1, \ldots, f_m)
    \Longrightarrow \gamma'(f_1, \ldots, f_m)\]
given by composing components of the associativity constraint $a$.
(Uniqueness follows from coherence for a bicategory.)

We then generate occupants of this niche as
    \[ \begin{array} {r} \{\alpha \circ u_{\gamma} : \alpha \in
    X(2)_1, s(\alpha) = \gamma(f_1, \ldots, f_m)\} \\ \  \end{array}
    \left/ \begin{array}{l} \\ \sim \end{array} \right.  \]
where $\sim$ is the equivalence relation generated by
\begin{enumerate}
\item $\alpha \circ u_{\gamma} = \beta \circ u_{\gamma'} \iff
\beta \circ a_{\gamma \gamma'} = \alpha \ \in \cl{B}$
\item $1 \circ u_\gamma = u_\gamma$.
\end{enumerate}

\noindent Note in particular that since $1 \circ a_{\gamma\gamma'}
= a_{\gamma\gamma'}$ we have
    \[a_{\gamma\gamma'}\circ u_{\gamma} = u_{\gamma'}.\]
So, given any $\gamma$, every occupant of the niche is uniquely
expressible as $\alpha \circ u_{\gamma}$, with $\alpha \in
X(2)_1$.  This shows that $u_\gamma$ is indeed universal, and
completes the definition of $X(2)$.

Note that the universality of the $u_\gamma$ follows from
coherence for classical bicategories, as it depends on the fact
that any two composites of given 1-cells are uniquely isomorphic.

We now construct $X(3)$.  We must specify a unique 3-cell for any
3-niche, that is, a unique composite 2-cell for any formal
composite of 2-cells.

\renewcommand{\labelenumi}{\arabic{enumi})}
\begin{enumerate}
\item First, composites of 1-ary 2-cells are determined by 2-cell
composition in \cl{B}.
\item Next we consider any composite of the form $c \circ \iota$.
We define the composites by

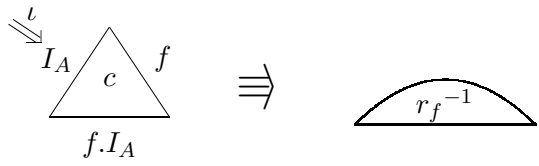
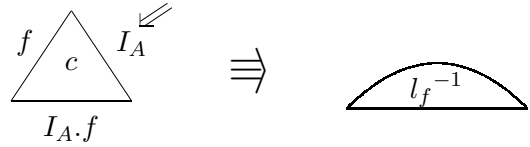
\begin{center}
\begin{picture}(30,22)(-9,-4)
\diagbw{$$}{$$}{$$}{$I_A$}{$f$}{$f.I_A$}{$\iota$}{$c$}
\end{picture}
\begin{picture}(12,18) \diagp{5}{8}{6mm} \end{picture}
\begin{picture}(35,15)(5,-3)
\onetwob{10}{0}{1mm}{$$}{$$}{$$}{$$}{${r_f}^{-1}$}
\end{picture}
\end{center}

\noindent and similarly

\begin{center}
\begin{picture}(30,22)(-5,-4)
\diagbx{$$}{$$}{$$}{$f$}{$I_A$}{$I_A.f$}{$$}{$c$}
\end{picture}
\begin{picture}(12,18) \diagp{5}{8}{6mm} \end{picture}
\begin{picture}(35,15)(5,-3)
\onetwob{10}{0}{1mm}{$$}{$$}{$$}{$$}{${l_f}^{-1}$}
\end{picture}.
\end{center}

\item Now consider a composite of the form

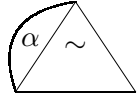
\begin{center}
\begin{picture}(25,15)(-5,-1)
\diagbv{$$}{$$}{$$}{$$}{$$}{$$}{$$}{$\alpha$}{$\sim$}
\end{picture}
\end{center}

\noindent where $\alpha$ is any 1-ary 2-cell.  We put

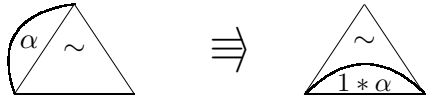
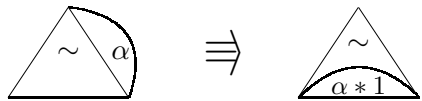
\begin{center}
\begin{picture}(25,15)(-5,-1)
\diagbv{$$}{$$}{$$}{$$}{$$}{$$}{$$}{$\alpha$}{$\sim$}
\end{picture}
\begin{picture}(11,15) \three{7}{7}{1mm}{} \end{picture}
\begin{picture}(25,15)(-5,-1)
\diagbr{$$}{$$}{$$}{$$}{$$}{$$}{$$}{$\sim$}{{\small $1\ast
\alpha$}}
\end{picture}
\end{center}

\noindent and similarly

\begin{center}
\begin{picture}(25,15)(-5,-1)
\diagbt{$$}{$$}{$$}{$$}{$$}{$$}{$$}{$\alpha$}{$\sim$}
\end{picture}
\begin{picture}(11,15) \three{7}{7}{1mm}{} \end{picture}
\begin{picture}(25,15)(-5,-1)
\diagbr{$$}{$$}{$$}{$$}{$$}{$$}{$$}{$\sim$}{{\small $\alpha \ast
1$}}
\end{picture}
\end{center}

\item Now consider a formal composite of chosen 2-ary 2-cells
$c_{gf}$.  Such a diagram uniquely determines a composite $\gamma$
in \cl{B} of its boundary 1-cells.  So we set the composite 2-cell
in $X$ to be $u_{\gamma}$.  Conversely, any 2-cell $u_{\gamma}$
thus arises as the composite of some 2-cells $c$.

\item Finally, since we require that 2-cell composition be
strictly associative, we have determined all 3-cells in $X$.  For,
using the above cases, any nullary, 2-ary or $m$-ary composite can
be written in the form

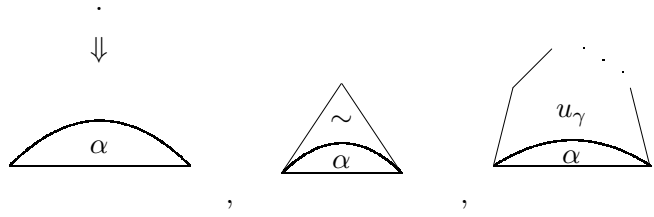
\begin{center}
\begin{picture}(25,30)(0,-5)
\diagbq{$$}{$$}{$$}{$$}{$\alpha$}
\end{picture}
\ \ ,
\begin{picture}(25,22)(-5,-4)
\diagbr{$$}{$$}{$$}{$$}{$$}{$$}{$$}{$\sim$}{$\alpha$}
\end{picture}
\ \ ,
\begin{picture}(25,25)(-2,-5)
\mtwob{$$}{$$}{$$}{$$}{$$} \diagcc{$u_\gamma$}{}{}{$\alpha$}
\end{picture}

\end{center}

\noindent respectively, where $\alpha$ is a composite of 1-ary
2-cells which we can then compose in \cl{B}.
\end{enumerate}


This completes the definition of the opetopic set $X$; it remains
to check that $X$ is 2-coherent. Certainly, every 3-niche has a
unique occupant by construction.  A 2-cell $\alpha \circ \iota$,
$\alpha \circ c$ or $\alpha \circ u_{\gamma}$ is universal if and
only if $\alpha$ is universal, that is, if and only if $\alpha$ is
invertible in \cl{B}. So every 2-niche has a universal occupant
and composites of universal 2-cells are universal.

We can check that a 1-cell in $X$ is universal if and only if it
is an (internal) equivalence in \cl{B}; this follows by an
analogous argument to the `Yoneda' proof of
Proposition~\ref{bicatyon}.  So every 1-niche has a universal
occupant $I_A$, and composites of universal 1-cells are universal.

So $X$ is a biased opetopic 2-category, with chosen universal
2-cells $\iota$ and $c$, and it is clear that $\zeta(X) = \cl{B}$.
So $\zeta$ is surjective.

\begin{itemize} \item We show that $\zeta$ is full.  \end{itemize}

Let $X$ and $X'$ be biased opetopic 2-categories, and suppose we
have a morphism of bicategories
    \[(G,\phi) : \cl{B}_X \lra \cl{B}_{X'}.\]
We define a morphism $F:X \lra X'$ as follows.  For $k=0$ and
$k=1$ the functions
    \[F_k:X(k)\lra X'(k)\]
are given by the action of $G$ on objects and 1-cells
respectively.  We construct $F_2$ as follows.  The action of $F_2$
on 1-ary 2-cells is the action of $G$ on 2-cells of $\cl{B}_X$.
For 0-ary 2-cells, we observe that any such is expressible
uniquely as

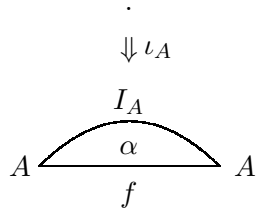
\begin{center}
\begin{picture}(25,30)(0,-5)
\diagbq{$A$}{$I_A$}{$f$}{$\iota_A$}{$\alpha$}
\end{picture}
\end{center}

\noindent where $\iota_A$ is the chosen universal for$X$.  Then we
define

\begin{center}
\begin{picture}(10,18) \put(0,7){$F_2 :$} \end{picture}
\begin{picture}(35,30)(-8,-5)
\diagbq{$A$}{$$}{$$}{$\iota_A$}{$\alpha$}
\end{picture}
\begin{picture}(25,18) \put(7,7){$\longmapsto$} \end{picture}
\begin{picture}(35,30)(0,-5)
\diagbq{$FA$}{$$}{$$}{$F\iota_A$}{$F\alpha$}
\end{picture}
\end{center}

\begin{picture}(12,18) \put(65,17){=} \end{picture}
\begin{picture}(45,40)(-70,-17)
\diagbq{$FA$}{$I_{FA}$}{}{$\iota_{FA}$}{$\phi_A$}
\qbezier(0,0)(12,-12)(24,0) %
\put(12,-7){\makebox(0,0)[t]{$Ff$}} %
\put(12,-2.5){\makebox(0,0)[c]{$F\alpha$}}
\end{picture}

\noindent where $\iota_{FA}$ is the appropriate chosen universal
for $X'$; this assignation is well-defined by uniqueness.

For $m\geq 2$, any $m$-ary 2-cell is expressible in the form

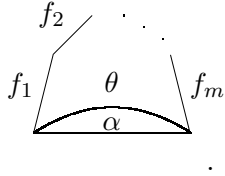
\begin{center}
\begin{picture}(25,25)(-2,-5)
\mtwob{$f_1$}{$f_2$}{$f_m$}{$$}{$$}
\diagcc{$\theta$}{}{}{$\alpha$}
\end{picture}.
\end{center}

\noindent Here $\theta$ is the composite of some configuration of
chosen universals $c$, determining a 1-cell composite $\gamma(f_1,
\ldots, f_m)$ in \cl{B}, and $\alpha : \gamma \Longrightarrow g$.
Then we define

\begin{center}
\begin{picture}(10,18) \put(0,7){$F_m :$} \end{picture}
\begin{picture}(25,25)(-2,-4)
\mtwob{}{}{}{}{} \diagcc{{\tiny $(c_1,c_2,\ldots)$}}{}{}{$\alpha$}
\end{picture}
\begin{picture}(18,18) \put(7,7){$\longmapsto$} \end{picture}
\begin{picture}(25,25)(-2,-4)
\mtwob{}{}{}{}{} \diagcc{{\tiny $(Fc_1,Fc_2,\ldots)$}}{}{}{{\tiny
$F\alpha$}}
\end{picture}
\end{center}

\begin{picture}(12,18) \put(65,17){=} \end{picture}
\begin{picture}(25,30)(-65,-10)
\mtwob{}{}{}{}{} \diagcc{{\tiny
$({c_1}',{c_2}',\ldots)$}}{}{}{{\tiny $\Phi$}} \diagcd{{\tiny
$F\alpha$}}
\end{picture}

\noindent where $\Phi$ is the appropriate composite of components
of the constraint $\phi$.  This assignation is well-defined by
uniqueness and the axioms for a morphism of bicategories.

It is clear from the construction that this is a morphism of
biased opetopic 2-categories, and that
    \[\zeta(F) = (G,\phi).\]
So $\zeta$ is full.

\begin{itemize} \item We show that $\zeta$ is faithful.
\end{itemize}

Consider morphisms $F, F'$ of unbiased opetopic 2-categories, such
that $\zeta(F)=\zeta(F')$.  Write $\zeta(F)=(G,\phi)$ and
$\zeta(F')=(G',\phi')$.

Certainly since $G=G'$ on objects and 1-cells we have $F_0=F'_0$
and $F_1=F'_1$. Similarly, $G=G'$ on (bicategorical) 2-cells gives
$F_2=F'_2$ on (opetopic) 1-ary 2-cells.  For $m$-ary 2-cells with
$m \neq 1$ consider again the above presentation of 2-cells.  Then
$\phi=\phi'$ gives $F_2 = F'_2$ on all opetopic 2-cells.  So
$\zeta$ is faithful.

\bigskip

So finally we may conclude that $\zeta$ exhibits an equivalence
    \[\cat{Opic-$2$-Cat}_b \simarr \cat{Bicat}\]
as required.
\end{prfof}

\begin{prfof}{Theorem \ref{bicatthmd1}}  By Proposition~\ref{bicatprope1} we
have
    \[\cat{Opic-$2$-Cat}_b \simarr \cat{Bicat}\]
and by Proposition~\ref{bicatpropf1} we have
    \[\cat{Opic-$2$-Cat}_b \simarr \cat{Opic-$2$-Cat}\]
so we have an equivalence
    \[\cat{Opic-$2$-Cat} \simeq \cat{Bicat}\]
as required.
\end{prfof}

\begin{remarks} \end{remarks}
\begin{enumerate}
\item Note that the final equivalence is not surjective in {\em
either} direction.  Left-to-right involves a choice of universal
2-cells; right-to-left involves generating sets of 3-cells and
$k$-ary 2-cells (for $k\neq 1$) which are only defined up to
isomorphism. Observe that a different choice of universal 2-cells
yields a bicategory non-trivially isomorphic but with the same
cells.

\item The term `biased opetopic 2-category' is used in the spirit of
Leinster's work on biased and unbiased bicategories (\cite{lei5}).
Rather than pick universal $m$-ary 2-cells for just $m=0,2$, we
might pick universals for all $m \geq 0$.  Again with no further
stipulations on morphisms, this yields an equivalent category of
`unbiased opetopic 2-categories'.  By a straightforward
modification of the above proof, we may see that this corresponds
to the theory of unbiased bicategories; Leinster has shown
directly that the biased and unbiased theories are equivalent.

\item In fact, we may choose any number of universal $m$-ary 2-cells for
each $m$ and define a category obviously equivalent to
\cat{Opic-$2$-Cat}, by making no stipulation on morphisms.  We
might then ask: when does this yield a theory of bicategories?  In
order to modify the above construction as required, we need enough
chosen universals to give a complete presentation of the 2-cells
of $X$.  From the observations in Section~\ref{bicatsecd1} we see
that this is possible provided we have chosen at least one 0-ary
universal, and at least one $m$-ary universal for some $m > 1$
(for each appropriate niche). This idea is discussed in
\cite{lei5} (Appendix A); in the opetopic setting it is immediate
that each resulting category of `bicategories' is equivalent.

\item Like Leinster, we might observe that the equivalence of {\em
categories}
    \[\cat{Opic-$2$-Cat} \simeq \cat{Bicat}\]
is two levels `better' than we might have asked; we have a
comparison at the 1-dimensional level without having to invoke 3-
or even 2-dimensional structures.  So the theory might already be
seen as fruitful despite the lack of an $(n+1)$-category of
$n$-categories.

\end{enumerate}

In summary, we have the following equivalences, surjective in the
directions shown:
     \[\cat{Opic-$2$-Cat} \stackrel{\sim}{\longleftarrow}
     \cat{Opic-$2$-Cat}_b \simarr \cat{Bicat}.\]


\subsection{Strictness} \label{bicatstr}

In this section we discuss (informally) various possible notions
of strictness in the opetopic setting, and compare these with the
classical biased and unbiased settings.

In the classical theory of bicategories, `strictness' (of
bicategories or their morphisms) is determined by the `strictness'
of the constraints; in general `lax' for plain morphisms, `weak'
for isomorphisms and `strict' for identities.

In the opetopic theory we cannot make such definitions, since we
do not have those constraints unless we have chosen universal
2-cells.   Even then the constraints are not explicitly given. So
we must define strictness by some other means; we may define
stricter and weaker notions in terms of universals.

We first turn our attention to morphisms.  Recall that the
original Baez-Dolan definition demanded that a morphism preserve
universality; this is stronger than the general morphisms we use
in our definition of \cat{Opic-$2$-Cat}.

\begin{proposition}\label{bicatpropg1} Recall (Proposition
\ref{bicatprope1}) that we have an equivalence
    \[\zeta:\cat{Opic-$2$-Cat}_b \simarr \cat{Bicat}.\]
Let $F$ be a morphism of opetopic 2-categories.  Then $F$
preserves universals iff $\zeta(F)$ is a weak functor
(homomorphism) of bicategories.
\end{proposition}

\begin{prf} Suppose $F:X\lra X'$ preserves universals.  Then the chosen
universal in $X$

\begin{center}
\begin{picture}(25,18)(-5,-1)
\twotwob{}{}{}{$f$}{$g$}{$gf$}{$c$}
\end{picture}
\end{center}
\noindent becomes, under the action of $F$, a universal in $X'$

\begin{center}
\begin{picture}(25,22)(-5,-4)
\twotwob{}{}{}{$Ff$}{$Fg$}{$F(gf)$}{$Fc$}
\end{picture}
\end{center}
\noindent inducing

\begin{center}
\begin{picture}(25,22)(-5,-4)
\diagbr{$$}{$$}{$$}{$Ff$}{$Fg$}{$$}{$$}{$Fc$}{{\small
$\phi^{-1}$}}
\end{picture}
\begin{picture}(12,18) \put(5,7){=} \end{picture}
\begin{picture}(25,22)(-5,-4)
\twotwob{$$}{$$}{$$}{$Ff$}{$Fg$}{$Fg.Ff$}{$\sim$}
\end{picture}
\end{center}

\noindent so $\phi_{ABC}$ is an isomorphism.

Conversely suppose $\phi_{gf}$ and $\phi_A$ are invertible for all
$f,g,A$.  First note that  1-ary universal 2-cells are always
preserved (clear from the case $n=1$).  Now, any universal can be
expressed as

\begin{center}
\begin{picture}(25,25)(-2,-5)
\mtwob{$$}{$$}{$$}{$$}{$$} \diagcc{$\theta$}{}{}{{\small
$\alpha$}}
\end{picture}
\end{center}

\noindent where $\theta$ is some composite of chosen universals
and $\alpha$ is universal.  Now applying F we have

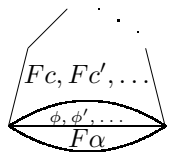
\begin{center}
\begin{picture}(25,25)(-2,-5)
\mtwob{$$}{$$}{$$}{$$}{$$} \diagcc{{\small
$Fc,Fc',\ldots$}}{}{}{{\tiny $\phi,\phi',\ldots$}} \diagcd{{\small
$F\alpha$}}
\end{picture}
\end{center}

\noindent which is universal since $F \alpha$ is universal.

The result for 1-cells follows (with some effort). \end{prf}

\begin{definition}
We write  \cat{Opic-$2$-Cat(weak)},
$\cat{Opic-$2$-Cat}_b$\cat{(weak)} and \cat{Bicat(weak)} for the
lluf subcategories with only weak morphisms.
\end{definition}

\begin{proposition} \label{bicatproph1} The equivalences given in the
proofs of Propositions~\ref{bicatpropf1} and \ref{bicatprope1}
restrict to equivalences
    \[\cat{Opic-$2$-Cat(weak)}
    \stackrel{\sim}{\longleftarrow}
    \cat{Opic-$2$-Cat}_b\cat{(weak)}
\simarr \cat{Bicat(weak)}\] surjective in the directions shown.
\end{proposition}

\begin{prf} The first equivalence is clear from the definitions
and the second follows from Proposition~\ref{bicatpropg1}.  Since
these are lluf subcategories the functors are clearly still
surjective.
\end{prf}

Since we have still made no stipulation about the action of
morphisms on chosen universals, it is clear that we will still
have a result of the form `all theories are equivalent' (cf
\cite{lei6}). That is, regardless of the number of universals
chosen, the category-with-weak-morphisms will remain equivalent to
the category \cat{Opic-$2$-Cat(weak)}.  This ceases to be so in
the strict case.

There is no obvious way of further strengthening the conditions
imposed on morphisms in \cat{Opic-$2$-Cat(weak)}, but if we
consider $\cat{Opic-$2$-Cat}_b\cat{(weak)}$, we can further demand
that chosen universals be preserved.


\begin{proposition} \label{bicatpropi1}  Let $F$ be a weak morphism of
biased opetopic 2-categories.  Then $F$ preserves chosen
universals iff $\zeta(F)$ is strict.
\end{proposition}

\begin{prf}
`$\Rightarrow$' is clear from the definition of $\zeta$. Now for
{\em any} morphism $(F,\phi)$ of opetopic 2-categories we have

\begin{center}
\begin{picture}(25,18)(0,-5)
\nulltwo{0}{0}{1mm}{$$}{$$}{$F\iota_A$}
\end{picture}
\begin{picture}(12,18) \put(5,7){=} \end{picture}
\begin{picture}(25,30)(0,-5)
\diagbq{$$}{$$}{$$}{$\iota_{FA}$}{$\phi_A$}
\end{picture}
\end{center}

\noindent and

\begin{center}
\begin{picture}(25,22)(-5,-4)
\twotwob{$$}{$$}{$$}{$$}{$$}{$$}{$Fc_{gf}$}
\end{picture}
\begin{picture}(12,18) \put(5,7){=} \end{picture}
\begin{picture}(25,22)(-5,-4)
\diagbr{$$}{$$}{$$}{$$}{$$}{$$}{$$}{$c_{Fg.Ff}$}{$\phi_{gf}$}
\end{picture}
\end{center}

\noindent so clearly if $(F,\phi)$ is strict then $F$ preserves
chosen universals.  \end{prf}

\begin{definition}  We call a weak morphism of biased opetopic
2-categories {\em strict} if it preserves chosen universal
2-cells.
\end{definition}

Write $\cat{Opic-$2$-Cat}_b\cat{(str)}$ and \cat{Bicat(str)} for
the lluf subcategories with only strict morphisms.

\begin{proposition}\label{bicatpropj1} The previously defined equivalence
restricts to an equivalence
    \[\cat{Opic-$2$-Cat}_b\cat{(str)} \simarr \cat{Bicat(str)}\]
surjective in the direction shown.
\end{proposition}

\begin{prf}  Follows immediately from Proposition~\ref{bicatpropi1}
\end{prf}

We now consider the possibility of altering the structures of the
2-categories themselves.  Considering the structures used so far
as `weak', we might try to find either lax or strict opetopic
$n$-categories.

In the lax direction we might consider removing the condition that
universals compose to universals.  Observe that in the case $n=1$
we do not use this condition to prove
    \[\cat{Opic-$1$-Cat} \simeq \cat{Cat}\]
so a `lax opetopic 1-category' would be just the same as a weak
one, as we would hope.

However, for $n=2$ it is not clear that this `laxification'
produces a useful structure for the general or biased theories.
Consider instead the case in which $m$-ary universal 2-cells have
been chosen for every $m \geq 0$.  That is, we define an `unbiased
opetopic 2-category' to be one in which {\em every} 2-niche has a
chosen universal occupant.

If we now remove the condition that composites of universals be
universals, we have certain 2-cell `constraints' induced by the
chosen universals.  For example we have

\begin{center}
\begin{picture}(25,20)(-2,-4)
\assleft{$f$}{$g$}{$h$}{$$}{$(hg)f$} \assleftb{$\sim$}{$\sim$}
\end{picture}
\begin{picture}(12,18) \put(5,7){=} \end{picture}
\begin{picture}(25,20)(-2,-4)
\threetwob{$$}{$$}{$$}{$$}{$f$}{$g$}{$h$}{$$}{$$}
\diagcc{$\sim$}{}{}{$\gamma$}
\end{picture}
\end{center}

\noindent and thus an induced 2-cell
    \[\gamma:hgf \Rightarrow (hg)f.\]
This produces a structure something like a `lax unbiased
bicategory' in the sense of Leinster (\cite{lei5}) except that the
constraints $\gamma$ are acting in the opposite direction.

For strictness there is likewise no obvious way of imposing
stronger conditions on an opetopic 2-category. Once we have chosen
universals, we might demand that the chosen universals compose to
chosen universals, but this will certainly not be possible unless
we have chosen $m$-ary universals for all $m \geq 0$.  So once
again we find ourselves in the unbiased theory.

If we have one chosen universal for each 2-niche, the above
condition forces strict associativity and left and right unit
action.  So we have a 2-category; this is to be expected since
Leinster has already observed that unbiased 2-categories are in
one-to-one correspondence with 2-categories.  (There is a
possibility of more interesting structure if a niche has more than
one chosen universal.)

\subsubsection*{Remarks}

From this informal discussion we see that the theory of opetopic
2-categories neither laxifies nor strictifies particularly
naturally.  In the lax direction, this is perhaps consistent with
the fact that there is no very satisfactory lax version of
classical bicategories.  In the strict direction, this
demonstrates why we have found it hard to state a coherence
theorem of the form `every bicategory is biequivalent to a
2-category'; we simply do not know what a `strict opetopic
2-category' is.  (Note however that statements of the form `all
diagrams commute' are much less problematic.)

We have already observed that there are (at least) two possible
ways of removing the bias in a bicategory: we may choose $m$-ary
composites for no $m$, or all $m$.  It appears that, although the
former philosophy may be viewed as being more egalitarian towards
all universal cells, the latter provides more footholds for
exploring the theory.

\subsection{Conclusions}

We might regard the category of opetopic 2-categories (with {\em
no} choices made) as being the most general of all the theories discussed in
this work. However we will also observe that in describing the
2-cells, performing calculations or exploring the theory further,
it is often more practical to use some presentation of 2-cells,
that is, to make choices of universals either explicitly or
implicitly.

In the opetopic setting the choice of universals is `free' in the
sense that no axioms are required; all axioms are subsumed into
the conditions for $n$-coherence.  So in each separate case the
axioms do not have to be stated explicitly.

This was suggested in \cite{bd1} as one of the motivations for the
opetopic approach to $n$-categories, since as $n$ increases, the
axioms for an $n$-category increase in complexity with fiendish
rapidity.  This work demonstrates a sense in which this idea is
realised.

\appendix

\renewcommand{\qbeziermax}{150}

\section{Calculations for Section~\ref{bicatbicat}}
\label{calc}

In this appendix we perform the calculations deferred from
Section~\ref{bicatbicat}.  However, we first introduce some
shorthand to deal with some of the more unwieldy parts of the
algebra.

\subsection{Shorthand for calculations}

The following shorthand is used for calculations in an opetopic
2-category.

\renewcommand{\labelenumi}{\roman{enumi})}
\begin{enumerate}

\item Since 3-niche occupants are unique, we may omit the target
of a 3-cell without ambiguity.  We then write an equality to
indicate that the 3-cells have the same target.  For example we
might write

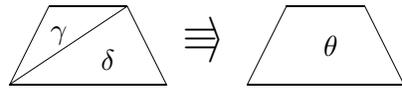
\begin{center}
\begin{picture}(25,20)(-2,-4)
\assleftb{$\beta$}{$\alpha$}
\end{picture}
\begin{picture}(12,18) \put(5,7){=} \end{picture}
\begin{picture}(25,20)(-2,-4)
\assrightb{$\gamma$}{$\delta$}
\end{picture}
\end{center}

\noindent meaning

\sunit
\begin{center}
\begin{picture}(25,20)(-2,-4)
\assleftb{$\beta$}{$\alpha$}
\end{picture}
\begin{picture}(4,18) \diagp{1}{10}{6mm} \end{picture}
\begin{picture}(25,20)(-2,-4)
\threetwob{}{}{}{}{}{}{}{}{$\theta$}
\end{picture} \end{center}

and

\begin{center}
\begin{picture}(25,20)(-2,-4)
\assrightb{$\gamma$}{$\delta$}
\end{picture}
\begin{picture}(4,18) \diagp{1}{10}{6mm} \end{picture}
\begin{picture}(25,20)(-2,-4)
\threetwob{}{}{}{}{}{}{}{}{$\theta$} \nolinebreak
\end{picture}.
\end{center}

\noindent

\item Recall that, by uniqueness of 3-niche occupants, we have
associativity of 2-cell composition.  So we may substitute `equal'
(in the above sense) 2-cell composites in part of the domain of
another 3-cell.  For example, given

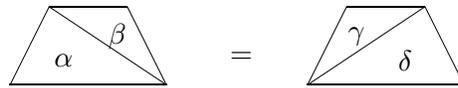
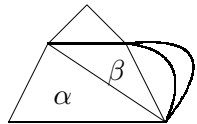
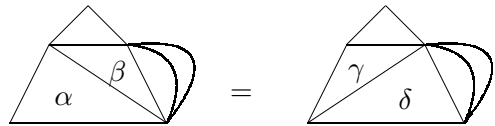
\begin{center}
\begin{picture}(25,20)(-2,-4)
\assleftb{$\beta$}{$\alpha$}
\end{picture}
\begin{picture}(12,18) \put(5,7){=} \end{picture}
\begin{picture}(25,20)(-2,-4)
\assrightb{$\gamma$}{$\delta$}
\end{picture}
\end{center}

\noindent and a 3-cell

\begin{center}
\begin{picture}(25,20)(-2,-4)
\assleftb{$\beta$}{$\alpha$} %
\diagcm %
\diagcg{} \diagch{}
\end{picture}
\end{center}

\noindent we have

\begin{center}
\begin{picture}(25,20)(-2,-4)
\assleftb{$\beta$}{$\alpha$} %
\diagcm %
\diagcg{} \diagch{}
\end{picture}
\begin{picture}(12,18) \put(5,7){=} \end{picture}
\begin{picture}(25,20)(-2,-4)
\assrightb{$\gamma$}{$\delta$} %
\diagcm %
\diagcg{} \diagch{}
\end{picture}.
\end{center}

\noindent This is shorthand for the following

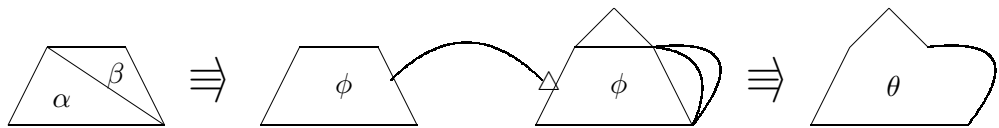
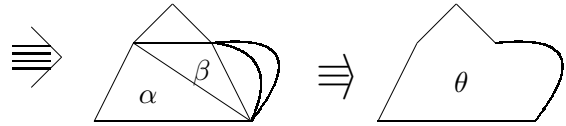
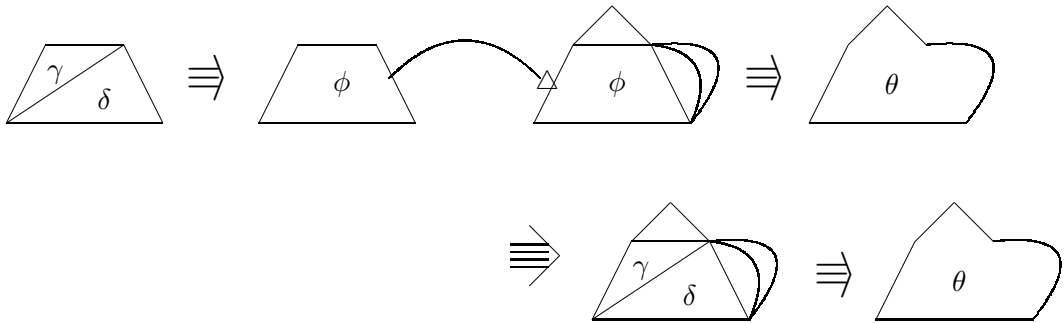
\begin{center}
\begin{picture}(5,25)(18,-4)
\assleftb{$\beta$}{$\alpha$}
\end{picture}
\begin{picture}(6,18) \three{2}{10}{1mm}{} \end{picture}
\begin{picture}(25,20)(-2,-4)
\threetwob{}{}{}{}{}{}{}{}{$\phi$}
\end{picture}
\begin{picture}(9,10)(2,0) \qbezier(-5,10)(5,20)(15,10)
\put(14.5,9.5){\makebox(0,0)[l]{{\small $\triangle$}}}
\end{picture}
\begin{picture}(25,20)(-2,-4)
\threetwob{}{}{}{}{}{}{}{}{$\phi$} \diagcm \diagcg{} \diagch{}
\end{picture}
\begin{picture}(9,18) \three{6}{10}{1mm}{} \end{picture}
\begin{picture}(25,20)(-2,-4)
\put(0,0){\setlength{\unitlength}{0.52mm}  %
\put(0,0){\line(1,0){40}} \put(0,0){\line(1,2){10}}} %
\put(11.05,5.2){\makebox(0,0)[c]{$\theta$}}
\diagcm{} \diagch{}  %
\end{picture}
\end{center} \nopagebreak

\begin{picture}(11,18) \diagc{50}{10}{1mm}{} \end{picture}
\begin{picture}(75,20)(-50,-4)
\assleftb{$\beta$}{$\alpha$} %
\diagcm %
\diagcg{} \diagch{}
\end{picture}
\begin{picture}(9,18) \three{6}{10}{1mm}{} \end{picture}
\begin{picture}(25,20)(-2,-4)
\put(0,0){\setlength{\unitlength}{0.52mm}  %
\put(0,0){\line(1,0){40}} \put(0,0){\line(1,2){10}}} %
\put(11.05,5.2){\makebox(0,0)[c]{$\theta$}}
\diagcm{} \diagch{}  %
\end{picture}

\noindent and

\begin{center}
\begin{picture}(5,25)(18,-4)
\assrightb{$\gamma$}{$\delta$}
\end{picture}
\begin{picture}(6,18) \three{2}{10}{1mm}{} \end{picture}
\begin{picture}(25,20)(-2,-4)
\threetwob{}{}{}{}{}{}{}{}{$\phi$}
\end{picture}
\begin{picture}(9,10)(2,0) \qbezier(-5,10)(5,20)(15,10)
\put(14.5,9.5){\makebox(0,0)[l]{{\small $\triangle$}}}
\end{picture}
\begin{picture}(25,20)(-2,-4)
\threetwob{}{}{}{}{}{}{}{}{$\phi$} \diagcm \diagcg{} \diagch{}
\end{picture}
\begin{picture}(9,18) \three{6}{10}{1mm}{} \end{picture}
\begin{picture}(25,20)(-2,-4)
\put(0,0){\setlength{\unitlength}{0.52mm}  %
\put(0,0){\line(1,0){40}} \put(0,0){\line(1,2){10}}} %
\put(11.05,5.2){\makebox(0,0)[c]{$\theta$}}
\diagcm{} \diagch{}  %
\end{picture} \nopagebreak
\end{center} \nopagebreak

\nopagebreak

\begin{picture}(11,18) \diagc{50}{10}{1mm}{} \end{picture}
\begin{picture}(75,20)(-50,-4)
\assrightb{$\gamma$}{$\delta$} %
\diagcm %
\diagcg{} \diagch{}
\end{picture}
\begin{picture}(9,18) \three{6}{10}{1mm}{} \end{picture}
\begin{picture}(25,20)(-2,-4)
\put(0,0){\setlength{\unitlength}{0.52mm}  %
\put(0,0){\line(1,0){40}} \put(0,0){\line(1,2){10}}} %
\put(11.05,5.2){\makebox(0,0)[c]{$\theta$}}
\diagcm{} \diagch{}  %
\end{picture}.

\noindent

\item Recall that 2-cell identities act as identities on $k$-ary
2-cells for all $k$ (not only 1-ary 2-cells), for example

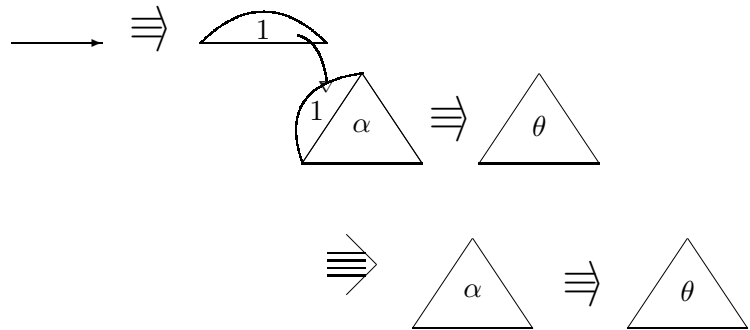
\begin{picture}(40,30)(-3,-20)
\nullthreeb{}{$1$}{}
\end{picture}
\begin{picture}(0,20) \diagad{5}{21}{1.5mm} \end{picture}
\begin{picture}(20,22)(-5,-4)
\diagbv{}{}{}{}{}{}{}{$1$}{} \twotwob{}{}{}{}{}{}{$\alpha$}
\end{picture}
\begin{picture}(6,18) \three{3}{10}{1mm}{} \end{picture}
\begin{picture}(15,22)(0,-4)
\twotwob{}{}{}{}{}{}{$\theta$}
\end{picture}

\begin{picture}(11,18) \diagc{50}{10}{1mm}{} \end{picture}
\begin{picture}(65,20)(-50,-4)
\twotwob{}{}{}{}{}{}{$\alpha$}
\end{picture}
\begin{picture}(9,18) \three{6}{10}{1mm}{} \end{picture}
\begin{picture}(25,20)(-2,-4)
\twotwob{}{}{}{}{}{}{$\theta$}
\end{picture}

\noindent so we have $\alpha=\theta$, that is

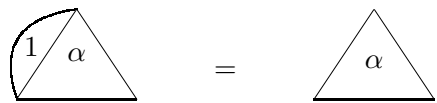
\begin{center}
\begin{picture}(25,22)(-5,-4)
\diagbv{$$}{$$}{$$}{$$}{$$}{$$}{$$}{$1$}{$\alpha$}
\end{picture}
\begin{picture}(12,18) \put(5,7){=} \end{picture}
\begin{picture}(25,22)(-5,-4)
\twotwob{$$}{$$}{$$}{$$}{$$}{$$}{$\alpha$}
\end{picture}.
\end{center}

\noindent

\item Note that if $u$ is any universal 2-cell, we have

\begin{center}
\begin{picture}(25,25)(-2,-5)
\mtwob{$$}{$$}{$$}{$$}{$$} \diagcc{$u$}{}{}{$\theta$}
\end{picture}
\begin{picture}(12,18) \put(5,7){=} \end{picture}
\begin{picture}(25,25)(-2,-5)
\mtwob{$$}{$$}{$$}{$$}{$$} \diagcc{$u$}{}{}{$\phi$}
\end{picture}
\begin{picture}(12,18) \put(5,7){$\Rightarrow\ \  \theta=\phi$} \end{picture}
\end{center}

\noindent by definition of universality.  This also holds if
$\theta$ and $\phi$ are 2-cell composites, for example

\begin{center}
\begin{picture}(16,10)(1,-4)
\diagco{1mm}{$\theta$} %
\end{picture}
\begin{picture}(6,10) \put(1,4){=} \end{picture}
\begin{picture}(16,10)(1,-4)
\diagco{1mm}{$\alpha$} %
\diagby{1mm}{$\beta$}%
\end{picture}
\begin{picture}(14,10) \put(5,4){$\Rightarrow$} \end{picture}
\begin{picture}(16,10)(1,-4)
\diagco{1mm}{$\alpha$} %
\diagby{1mm}{$\beta$}%
\end{picture}
\begin{picture}(6,10) \put(1,4){=} \end{picture}
\begin{picture}(16,10)(1,-4)
\diagco{1mm}{$\phi$} %
\end{picture}
\end{center}

\noindent and

\begin{center}
\begin{picture}(15,25)(8,-4)
\mtwob{}{}{}{}{} \diagcc{}{}{}{$\alpha$} \diagcd{$\beta$}
\end{picture}
\begin{picture}(6,10) \put(1,4){=} \end{picture}
\begin{picture}(25,25)(-2,-4)
\mtwob{}{}{}{}{} \diagcc{}{}{}{$\gamma$} \diagcd{$\delta$}
\end{picture}
\begin{picture}(14,10) \put(5,4){$\Rightarrow$} \end{picture}
\begin{picture}(16,10)(1,-4)
\diagco{1mm}{$\alpha$} %
\diagby{1mm}{$\beta$}%
\end{picture}
\begin{picture}(6,10) \put(1,4){=} \end{picture}
\begin{picture}(16,10)(1,-4)
\diagco{1mm}{$\gamma$} %
\diagby{1mm}{$\delta$}%
\end{picture}.
\end{center}

\noindent Furthermore, this holds if $u$ is a composite of
universals, since a composite of universals is universal, for
example if $u_1$ and $u_2$ are universal then

\begin{center}
\begin{picture}(25,20)(-2,-4)
\assrightb{$u_1$}{$$} \diagcc{}{}{$u_2$}{$\alpha$}
\end{picture}
\begin{picture}(12,18) \put(5,7){=} \end{picture}
\begin{picture}(25,20)(-2,-4)
\assrightb{$u_1$}{$$} \diagcc{}{}{$u_2$}{$\beta$}
\end{picture}
\begin{picture}(12,18) \put(5,7){$\Rightarrow\ \  \alpha=\beta$} \end{picture}
\end{center}

\noindent and in particular

\begin{center}
\begin{picture}(18,22)(0,-4)
\diagbr{$$}{$$}{$$}{$$}{$$}{$$}{$$}{$u$}{$\alpha$}
\end{picture}
\begin{picture}(2,18) \put(-1,7){=} \end{picture}
\begin{picture}(18,22)(0,-4)
\twotwob{$$}{$$}{$$}{$$}{$$}{$$}{$u$}
\end{picture}
\begin{picture}(8,18) \put(2,7){$\Rightarrow$} \end{picture}
\begin{picture}(20,15)(10,-3)
\onetwob{10}{0}{0.7mm}{$$}{$$}{$$}{$$}{$\alpha$}
\end{picture}
\begin{picture}(2,18) \put(-1,7){=} \end{picture}
\begin{picture}(2,15)(10,-3)
\onetwob{10}{0}{0.7mm}{$$}{$$}{$$}{$$}{$1$}
\end{picture}
\end{center}

\end{enumerate}


\subsection{Calculations}

Throughout this subsection, we use the notation and constructions
exactly as given in Section~\ref{bicatbicat}.

\begin{lemma} \label{bicatlemma1}
\begin{enumerate}
\item $1_g \ast 1_f = 1_{gf}$
\item $(\beta_2 \circ \beta_1) \ast (\alpha_2 \circ \alpha_1) =
(\beta_2 \ast \alpha_2) \circ (\beta_1 \ast \alpha_1)$ (middle 4
interchange)
\end{enumerate}
\end{lemma}

\begin{prf} \
\begin{enumerate}
\item We have

\begin{center}
\begin{picture}(25,22)(-5,-4)
\twotwob{$$}{$$}{$$}{$$}{$$}{$$}{$\sim$}
\end{picture}
\begin{picture}(12,18) \put(5,7){=} \end{picture}
\begin{picture}(25,22)(-5,-4)
\diagbt{$$}{$$}{$$}{$$}{$$}{$$}{$$}{$1$}{$\sim$}
\diagbv{$$}{$$}{$$}{$$}{$$}{$$}{$$}{$1$}{$$}
\end{picture}
\begin{picture}(12,18) \put(5,7){=} \end{picture}
\begin{picture}(25,22)(-5,-4)
\diagbr{$$}{$$}{$$}{$$}{$$}{$$}{$$}{$\sim$}{$1\ast 1$}
\end{picture}
\end{center}

\noindent by the action of 1 and definition of $\ast$, so

\begin{center}
\begin{picture}(35,15)(5,-3)
\onetwob{10}{0}{1mm}{$$}{$$}{$$}{$$}{$1\ast 1$}
\end{picture}
\begin{picture}(6,15) \put(4,7){=} \end{picture}
\begin{picture}(35,15)(5,-3)
\onetwob{10}{0}{1mm}{$$}{$$}{$$}{$$}{$1$}
\end{picture}
\end{center}

\noindent as required.

\item Given

\begin{center}
\begin{picture}(50,20)(0,-10)
\onetwoc{$f_1$}{$f_2$}{$f_3$}{$\alpha_1$}{$\alpha_2$} %
\put(24,0){\onetwoc{$g_1$}{$g_2$}{$g_3$}{$\beta_1$}{$\beta_2$}}
\end{picture}
\end{center}

\noindent we write

\begin{center}
\begin{picture}(25,22)(-5,-4)
\twotwob{}{}{}{$f_1$}{$g_1$}{$g_1f_1$}{$u_1$}
\end{picture}
\ ,
\begin{picture}(25,22)(-5,-4)
\twotwob{}{}{}{$f_2$}{$g_2$}{$g_2f_2$}{$u_2$}
\end{picture}
\ ,
\begin{picture}(25,22)(-5,-4)
\twotwob{}{}{}{$f_3$}{$g_3$}{$g_3f_3$}{$u_3$}
\end{picture}
\end{center}

\noindent for the chosen universal 2-cells as shown. Then we have

\begin{tiny}
\begin{center}
\begin{picture}(25,22)(-5,-4)
\diagca %
\diagcb %
\diagbv{$$}{$$}{$$}{$\alpha_1$}{$$}{$\beta_2$}{$$}{$\alpha_2$}{$u_3$} %
\diagbt{$$}{$$}{$$}{$$}{$\beta_1$}{$$}{$$}{$$}{$$} %
\end{picture}
\begin{picture}(12,18) \put(5,7){=} \end{picture}
\begin{picture}(25,22)(-5,-4)
\diagbv{$$}{$$}{$$}{$$}{$$}{$\beta_2\circ \beta_1$}{$$}{$\alpha_2\circ \alpha_1$}{$u_3$} %
\diagbt{$$}{$$}{$$}{$$}{$$}{$$}{$$}{$$}{$$} %
\end{picture}
\begin{picture}(12,18) \put(5,7){=} \end{picture}
\begin{picture}(25,22)(-5,-4)
\diagbr{$$}{$$}{$$}{$$}{$$}{$$}{$$}{$u_1$}{{\tiny $(\beta_2 \circ
\beta_1) \ast (\alpha_2 \circ \alpha_1)$}}
\end{picture}
\end{center}
\end{tiny}

\noindent by definition, but also

\begin{tiny}
\begin{center}
\begin{picture}(25,22)(-5,-4)
\diagca %
\diagcb %
\diagbv{$$}{$$}{$$}{$\alpha_1$}{$$}{$\beta_2$}{$$}{$\alpha_2$}{$u_3$} %
\diagbt{$$}{$$}{$$}{$$}{$\beta_1$}{$$}{$$}{$$}{$$} %
\end{picture}
\begin{picture}(12,18) \put(5,7){=} \end{picture}
\begin{picture}(25,22)(-5,-4)
\diagbv{$$}{$$}{$$}{$$}{$$}{$\beta_1$}{$$}{$\alpha_1$}{$$} %
\diagbt{$$}{$$}{$$}{$$}{$$}{$$}{$$}{$$}{$$} %
\diagbr{}{}{}{}{}{}{}{$u_2$}{$\beta_2 \ast \alpha_2$}
\end{picture}
\begin{picture}(12,18) \put(5,7){=} \end{picture}
\begin{picture}(25,22)(-5,-4)
\diagbr{$$}{$$}{$$}{$$}{$$}{$$}{$$}{$u_1$}{{\tiny $(\beta_2 \ast
\alpha_2) \circ (\beta_1 \ast \alpha_1)$}}
\end{picture}
\end{center}
\end{tiny}

\noindent by definition, hence the result.\end{enumerate}
\end{prf}

\begin{lemma}\label{bicatlemma2} $a$ is natural
\end{lemma}

\begin{prf} Given 2-cells

\begin{center}
\begin{picture}(82,15)(5,-3)
\onetwob{10}{0}{1mm}{$$}{$$}{$f_1$}{$f_2$}{$\alpha$}
\onetwob{34}{0}{1mm}{$$}{$$}{$g_1$}{$g_2$}{$\beta$}
\onetwob{58}{0}{1mm}{$$}{$$}{$h_1$}{$h_2$}{$\gamma$}
\end{picture}
\end{center}

\noindent we need to show that the following naturality square
commutes

\begin{center}
\setlength{\unitlength}{0.2em}
\begin{picture}(55,40)(2,5)    %

\put(10,10){\makebox(0,0){$(h_2 g_2)f_1$}}  
\put(10,35){\makebox(0,0){$(h_1 g_1)f_1$}}  
\put(45,35){\makebox(0,0){$h_1(g_1 f_1)$}}  
\put(45,10){\makebox(0,0){$h_2(g_2 f_2)$}}  

\put(19,35){\vector(1,0){15}}  
\put(19,10){\vector(1,0){15}}  
\put(10,30){\vector(0,-1){15}} 
\put(45,30){\vector(0,-1){15}} 

\put(8,23){\makebox(0,0)[r]{$(\gamma \ast \beta)\ast \alpha$}} 
\put(47,23){\makebox(0,0)[l]{$\gamma \ast (\beta \ast \alpha)$}} 
\put(27,37){\makebox(0,0)[b]{$a$}} 
\put(27,8){\makebox(0,0)[t]{$a$}} 

\put(60,9){\makebox(0,0){.}} 

\end{picture}
\end{center}

\sunit

We have

\begin{tiny}
\begin{center}
\begin{picture}(25,20)(-2,-4)
\assleftb{}{} \diagcf{$\alpha$} \diagci{$\beta$} \diagcg{$\gamma$}
\diagcd{$a$}
\end{picture}
\begin{picture}(12,18) \put(5,7){=} \end{picture}
\begin{picture}(25,20)(-2,-4)
\assleftb{}{} \diagcf{$\alpha$} \diagcd{$a$} \diagck{}{{\tiny
$\gamma \ast \beta$}}
\end{picture}
\begin{picture}(12,18) \put(5,7){=} \end{picture}
\begin{picture}(25,20)(-2,-4)
\assleftb{}{} \diagcd{$a$} \diagcc{}{}{}{$(\gamma \ast \beta) \ast
\alpha$}
\end{picture} \nopagebreak
\begin{picture}(25,30)(-2,-4) \put(10,20){$\parallel$}
\assrightb{}{} \diagcf{$\alpha$} \diagci{$\beta$}
\diagcg{$\gamma$}
\end{picture}
\begin{picture}(12,18) \put(5,7){=} \end{picture}
\begin{picture}(25,20)(-2,-4)
\assrightb{}{} \diagcg{$\gamma$} \diagcj{}{{\tiny $\beta \ast
\alpha$}}
\end{picture}
\begin{picture}(12,18) \put(5,7){=} \end{picture}
\begin{picture}(25,20)(-2,-4)
\assleftb{}{} \diagcd{$\gamma \ast (\beta \ast \alpha)$}
\diagcc{}{}{}{$a$}
\end{picture}

\end{center}
\end{tiny}

\noindent so by uniqueness we have

\begin{center}
\begin{picture}(37,22)(-7,-10)
\diagco{2mm}{$(\gamma \ast \beta) \ast \alpha$} %
\diagby{2mm}{$a$}%
\end{picture}
\begin{picture}(12,18) \put(10,11){=} \end{picture}
\begin{picture}(37,22)(-7,-10)
\diagco{2mm}{$a$} %
\diagby{2mm}{$\gamma \ast (\beta \ast \alpha)$}%
\end{picture}
\end{center}

\noindent as required. \end{prf}


\begin{lemma}\label{bicatlemma3} $r$ is natural
\end{lemma}

\begin{prf} Given a 2-cell

\begin{center}
\begin{picture}(35,15)(5,-3)
\onetwob{10}{0}{1mm}{$A$}{$B$}{$f_1$}{$f_2$}{$\alpha$}
\end{picture}
\end{center}

\noindent we need to show that the following naturality square
commutes

\begin{center}
\setlength{\unitlength}{0.2em}
\begin{picture}(55,40)(2,5)     %

\put(10,10){\makebox(0,0){$f_2 \circ I_A$}}  
\put(10,35){\makebox(0,0){$f_1 \circ I_A$}}  
\put(45,35){\makebox(0,0){$f_1$}}  
\put(45,10){\makebox(0,0){$f_2$}}  

\put(18,35){\vector(1,0){21}}  
\put(18,10){\vector(1,0){21}}  
\put(10,30){\vector(0,-1){15}} 
\put(45,30){\vector(0,-1){15}} 

\put(8,23){\makebox(0,0)[r]{$\alpha \ast 1$}} 
\put(47,23){\makebox(0,0)[l]{$\alpha$}} 
\put(27,37){\makebox(0,0)[b]{$r$}} 
\put(27,8){\makebox(0,0)[t]{$r$}} 

\put(55,9){\makebox(0,0){.}} 

\end{picture}
\end{center}
\sunit

\noindent Writing chosen composites as

\begin{center}
\begin{picture}(25,22)(-5,-4)
\twotwob{}{}{}{$I_A$}{$f_1$}{$f_1.I_A$}{$u_1$}
\end{picture}
\ ,
\begin{picture}(25,22)(-5,-4)
\twotwob{}{}{}{$I_A$}{$f_2$}{$f_2.I_A$}{$u_2$}
\end{picture}
\end{center}
\noindent we have

\begin{picture}(30,27)(-9,-4)
\diagbw{}{}{}{}{}{}{}{} \diagbr{}{}{}{}{}{}{}{$u_2$}{$r$}
\diagbt{}{}{}{}{}{}{}{$\alpha$}{}
\end{picture}
\begin{picture}(12,18) \put(5,7){=} \end{picture}
\begin{picture}(30,27)(-9,-4)
\diagbw{}{}{}{}{}{}{}{} \diagby{1mm}{$r$}
\diagbr{}{}{}{}{}{}{}{$u_1$}{$\alpha \ast 1$}
\end{picture}

\begin{picture}(29,15)(-9,-7) \put(7,9){$\parallel$}
\diagco{1mm}{$\alpha$} \diagby{1mm}{$1$}
\end{picture}
\begin{picture}(12,18) \put(6,7){=} \end{picture}
\begin{picture}(29,15)(-9,-7)
\diagco{1mm}{$1$} \diagby{1mm}{$\alpha$}%
\end{picture}
\begin{picture}(12,18) \put(5,7){=} \end{picture}
\begin{picture}(30,22)(-9,-4)
\diagbw{}{}{}{}{}{}{}{} \diagbr{}{}{}{}{}{}{}{$u_1$}{$r$}
\diagby{1mm}{$\alpha$}
\end{picture}

\noindent so by uniqueness

\begin{center}
\begin{picture}(25,15)(-5,-4)
\diagco{1mm}{$\alpha \ast 1$} %
\diagby{1mm}{$r$}%
\end{picture}
\begin{picture}(12,10) \put(5,7){=} \end{picture}
\begin{picture}(25,15)(-5,-4)
\diagco{1mm}{$r$} %
\diagby{1mm}{$\alpha$}%
\end{picture}
\end{center}

\noindent as required.
\end{prf}

\begin{lemma}\label{bicatlemma4}  $a$, $l$ and $r$ satisfy the axioms
for a bicategory.
\end{lemma}

\begin{prf} \
\begin{enumerate}
\item associativity pentagon

\setleng{0.52mm}
\begin{picture}(35,45)(10,-7)
\scalecp \put(40,0){\line(-2,3){20}} \put(40,0){\line(-2,1){40}}
\end{picture}
\begin{picture}(8,18) \put(0,17){=} \end{picture}
\begin{picture}(45,45)(0,-7)
\scalecp \put(40,0){\line(-2,3){20}} \put(0,0){\line(2,3){20}}
\diagcc{}{}{}{$a$}
\end{picture}
\begin{picture}(8,18) \put(0,17){=} \end{picture}
\begin{picture}(45,45)(0,-7)
\scalecp \put(0,0){\line(2,3){20}} \put(0,0){\line(2,1){40}}
\diagcc{}{}{}{$a$} \diagcd{$a$}
\end{picture}

\begin{picture}(35,55)(10,-7) \put(20,40){$\parallel$}
\scalecp \put(0,20){\line(1,0){40}} \put(40,0){\line(-2,1){40}}
\qbezier(0,20)(25,20)(40,0) \put(21,10){$a$}
\end{picture}
\begin{picture}(8,18) \put(0,17){=} \end{picture}
\begin{picture}(45,45)(0,-7)
\scalecp \put(0,20){\line(1,0){40}} \put(40,0){\line(-2,1){40}}
\diagcc{}{}{}{$a \ast 1$}
\end{picture}
\begin{picture}(8,18) \put(0,17){=} \end{picture}
\begin{picture}(45,45)(0,-7)
\scalecp \put(0,20){\line(1,0){40}} \put(0,0){\line(2,1){40}}
\diagcc{}{}{}{$a$} \diagcd{$a \ast 1$}
\end{picture}
\begin{picture}(8,18) \put(0,17){=} \end{picture}
\begin{picture}(50,45)(0,-7)
\scalecp \put(0,0){\line(2,3){20}} \put(0,0){\line(2,1){40}}
\qbezier(0,0)(15,20)(40,20) \put(19,11){$a$} \diagcc{}{}{}{$a$}
\diagcd{$a \ast 1$}
\end{picture}

\begin{picture}(165,65) \put(154,37){=} \end{picture}
\begin{picture}(45,65)(2,-27)
\scalecp \put(0,0){\line(2,3){20}} \put(0,0){\line(2,1){40}}
\diagcc{}{}{}{$1 \ast \alpha$} \diagcd{$a$} \diagce{$a \ast 1$}
\end{picture}

\sunit

\noindent so

\begin{center}
\begin{picture}(25,20)(-5,-9)
\diagco{1mm}{$a$} %
\diagby{1mm}{$a$}%
\put(8,7){\makebox(0,0)[c]{$((kh)g)f$}}
\put(8,-7){\makebox(0,0)[c]{$k(h(gf))$}}
\end{picture}
\begin{picture}(12,10) \put(5,9){=} \end{picture}
\begin{picture}(25,20)(-5,-9)
\diagco{1mm}{$1 \ast a$} %
\diagby{1mm}{$a$}%
\diagbz{$a \ast 1$} %
\end{picture}
\end{center}

\noindent as required.

\item unit triangle

\begin{picture}(20,25)(3,-4)
\assleftb{}{} \diagck{}{$r$} \diagcl{}
\end{picture}
\begin{picture}(6,18) \put(1,7){=} \end{picture}
\begin{picture}(20,20)(-2,-4)
\assleftb{}{} \diagcl{} \diagcc{}{}{}{$r \ast 1$}
\end{picture}

\begin{picture}(15,22)(2,-4) \put(8,16){$\parallel$}
\diagbt{}{}{}{}{}{}{}{1}{}
\end{picture}
\begin{picture}(15,18) \put(6,7){=} \end{picture}
\begin{picture}(21,22)(-2,-4)
\diagbv{}{}{}{}{}{}{}{1}{}
\end{picture}
\begin{picture}(6,18) \put(1,7){=} \end{picture}
\begin{picture}(22,20)(0,-4)
\assrightb{}{} \diagcl{} \diagcj{}{$l$}
\end{picture}
\begin{picture}(6,18) \put(1,7){=} \end{picture}
\begin{picture}(20,20)(0,-4)
\assrightb{}{} \diagcl{} \diagcc{}{}{}{$1 \ast l$}
\end{picture}

\begin{picture}(12,18) \put(86,7){=} \end{picture}
\begin{picture}(25,30)(-80,-4)
\assleftb{}{} \diagcl{} \diagcc{}{}{}{$a$} \diagcd{$1 \ast l$}
\end{picture}

\noindent so

\begin{center}
\begin{picture}(25,15)(-5,-4)
\diagco{1mm}{$r\ast 1$}
\end{picture}
\begin{picture}(12,18) \put(5,7){=} \end{picture}
\begin{picture}(25,15)(-5,-4)
\diagco{1mm}{$a$} %
\diagby{1mm}{$1 \ast l$}%
\end{picture}
\end{center}

\noindent as required.
\end{enumerate}
\end{prf}

\begin{lemma}\label{bicatlemma5} $\phi$ is natural.
\end{lemma}


\begin{prf} Given 2-cells

\begin{center}
\begin{picture}(58,15)(5,-3)
\onetwob{10}{0}{1mm}{$$}{$$}{$f_1$}{$f_2$}{$\alpha$}
\onetwob{34}{0}{1mm}{$$}{$$}{$g_1$}{$g_2$}{$\beta$}
\end{picture}
\end{center}

\noindent we need to show that the following diagram commutes

\begin{center}
\setlength{\unitlength}{0.2em}
\begin{picture}(65,40)(2,5)      %

\put(10,10){\makebox(0,0){$Fg_2 \circ Ff_2$}}  
\put(10,35){\makebox(0,0){$Fg_1 \circ Ff_1$}}  
\put(57,35){\makebox(0,0){$F(g_1 \circ f_1)$}}  
\put(57,10){\makebox(0,0){$F(g_2 \circ f_2)$}}  

\put(22,35){\vector(1,0){21}}  
\put(22,10){\vector(1,0){21}}  
\put(10,30){\vector(0,-1){15}} 
\put(57,30){\vector(0,-1){15}} 

\put(8,23){\makebox(0,0)[r]{$F\beta \ast F\alpha$}} 
\put(59,23){\makebox(0,0)[l]{$F(\beta \ast \alpha)$}} 
\put(33,37){\makebox(0,0)[b]{$\phi_{g_1 f_1}$}} 
\put(33,8){\makebox(0,0)[t]{$\phi_{g_2 f_2}$}} 

\put(79,9){\makebox(0,0){.}} 

\end{picture}
\end{center}
\sunit

\noindent We write the chosen universal 2-cells as

\begin{center}
\begin{picture}(25,22)(-5,-4)
\twotwob{}{}{}{$f_1$}{$g_1$}{$g_1.f_1$}{$v_1$}
\end{picture}
\ ,
\begin{picture}(25,22)(-5,-4)
\twotwob{}{}{}{$f_2$}{$g_2$}{$g_2.f_2$}{$v_2$}
\end{picture}
\ ,
\begin{picture}(25,22)(-5,-4)
\twotwob{}{}{}{$Ff_1$}{$Fg_1$}{$Fg_1.Ff_1$}{$u_1$}
\end{picture}
\ ,
\begin{picture}(25,22)(-5,-4)
\twotwob{}{}{}{$Ff_2$}{$Fg_2$}{$Fg_2.Ff_2$}{$u_2$}
\end{picture}
\end{center}
\noindent so

\begin{center}
\begin{picture}(18,22)(0,-4)
\diagbr{$$}{$$}{$$}{$$}{$$}{$$}{$$}{$u_1$}{$\phi$}
\end{picture}
\begin{picture}(2,18) \put(-1,7){=} \end{picture}
\begin{picture}(18,22)(0,-4)
\twotwob{$$}{$$}{$$}{$$}{$$}{$$}{$Fv_1$}
\end{picture}
\ ,
\begin{picture}(18,22)(0,-4)
\diagbr{$$}{$$}{$$}{$$}{$$}{$$}{$$}{$u_2$}{$\phi$}
\end{picture}
\begin{picture}(2,18) \put(-1,7){=} \end{picture}
\begin{picture}(18,22)(0,-4)
\twotwob{$$}{$$}{$$}{$$}{$$}{$$}{$Fv_2$}
\end{picture}.
\end{center}

\noindent We have

\setlength{\leng}{1.5mm}
\begin{center}
\begin{picture}(33,32)(-7,-6)
\scalebv{}{}{}{}{}{}{}{$F\alpha$}{$u_2$}
\scalebt{}{}{}{}{}{}{}{$F\beta$}{} \scaleby{$\phi_{g_2f_2}$}
\end{picture}
\begin{picture}(12,18) \put(5,12){=} \end{picture}
\begin{picture}(33,32)(-7,-6)
\scalebr{}{}{}{}{}{}{}{$u_1$}{$F\beta \ast F\alpha$}
\scaleby{$\phi_{g_2f_2}$}
\end{picture}
\end{center}

\noindent in $X'$, and in $X$ we have

\begin{center}
\begin{picture}(33,32)(-7,-6)
\scalebv{}{}{}{}{}{}{}{$\alpha$}{$v_2$}
\scalebt{}{}{}{}{}{}{}{$\beta$}{}
\end{picture}
\begin{picture}(12,18) \put(5,12){=} \end{picture}
\begin{picture}(33,32)(-7,-6)
\scalebr{}{}{}{}{}{}{}{$v_1$}{$\beta \ast \alpha$}
\end{picture}
\end{center}
\noindent so applying $F$, we have, since $F$ is strictly
functorial on 2-cells,

\begin{center}
\begin{picture}(34,32)(-7,-6)
\scalebv{}{}{}{}{}{}{}{$F\alpha$}{$Fv_2$}
\scalebt{}{}{}{}{}{}{}{$F\beta$}{}
\end{picture}
\begin{picture}(5,18) \put(4,12){=} \end{picture}
\begin{picture}(34,32)(-7,-6)
\scalebr{}{}{}{}{}{}{}{$Fv_1$}{$F(\beta \ast \alpha)$}
\end{picture}
\begin{picture}(5,18) \put(2,12){=} \end{picture}
\begin{picture}(37,32)(-7,-6)
\scalebr{}{}{}{}{}{}{}{$u_1$}{$\phi_{g_1f_1}$} \scaleby{$F(\beta
\ast \alpha)$}
\end{picture}
\end{center}

\begin{picture}(29,32)(-2,-6) \put(12,25){$\parallel$}
\scalebv{}{}{}{}{}{}{}{$F\alpha$}{$u_2$}
\scalebt{}{}{}{}{}{}{}{$F\beta$}{} \scaleby{$\phi_{g_2f_2}$}
\end{picture}
\begin{picture}(5,18) \put(4,12){=} \end{picture}
\begin{picture}(34,32)(-7,-6)
\scalebr{}{}{}{}{}{}{}{$u_1$}{$F\beta \ast F\alpha$}
\scaleby{$\phi_{g_2f_2}$}
\end{picture}

\noindent so by uniqueness we have

\setlength{\leng}{2mm} \setlength{\unitlength}{\leng}
\begin{center}
\begin{picture}(20,10)(-2.5,-4)
\scaleco{$F\beta \ast F\alpha$} %
\scaleby{$\phi_{g_2f_2}$}%
\end{picture}
\begin{picture}(6,10) \put(3,5){=} \end{picture}
\begin{picture}(20,10)(-2.5,-4)
\scaleco{$\phi_{g_1f_1}$} %
\scaleby{$F(\beta \ast \alpha)$}%
\end{picture}
\end{center}

\noindent as required. \end{prf}

\begin{lemma}\label{bicatlemma6} $(F,\phi)$ satisfies the axioms for a
morphism of bicategories.
\end{lemma}

\begin{prf} We have in $X$

\sunit
\begin{center}
\begin{picture}(25,20)(-2,-4)
\assleft{$f$}{$g$}{$h$}{}{} \assleftb{$\sim$}{$\sim$} \diagcd{$a$}
\end{picture}
\begin{picture}(12,18) \put(5,7){=} \end{picture}
\begin{picture}(25,20)(-2,-4)
\assrightb{$\sim$}{$\sim$}
\end{picture}
\end{center}

\noindent so applying $F$, we get in $X'$

\sunit
\begin{small}
\begin{center}
\begin{picture}(25,20)(-2,-4)
\assleft{$Ff$}{$Fg$}{$Fh$}{}{} \assleftb{}{} \diagcd{$Fa$}
\diagcc{}{$\sim$}{}{$\phi$} \diagck{$\sim$}{$\phi$}
\end{picture}
\begin{picture}(12,18) \put(5,7){=} \end{picture}
\begin{picture}(25,20)(-2,-4)
\assrightb{}{} \diagcj{$\sim$}{{\small $\phi$}}
\diagcc{}{}{$\sim$}{$\phi$}
\end{picture} \nopagebreak

\begin{picture}(25,30)(-2,-8) \put(10,16){$\parallel$}
\assleftb{$\sim$}{} \diagcd{$\phi$} \diagcc{}{$\sim$}{}{$\phi \ast
1$} \diagce{$Fa$}
\end{picture}
\begin{picture}(12,18) \end{picture}
\begin{picture}(25,30)(-2,-8) \put(10,16){$\parallel$}
\assleftb{$\sim$}{} \diagcd{$1 \ast \phi$}
\diagcc{}{$\sim$}{}{$a$} \diagce{$\phi$}
\end{picture}
\end{center}
\end{small}

\noindent as required. For $r$ we have in $X$

\begin{center}
\begin{picture}(25,22)(-5,-4)
\diagbr{}{}{}{}{}{}{}{$\sim$}{$r$} \diagbw{}{}{}{}{}{}{$\iota$}{}
\end{picture}
\begin{picture}(12,18) \put(5,7){=} \end{picture}
\begin{picture}(25,22)(-5,-4)
\diagco{1mm}{$1$}
\end{picture}
\end{center}

\noindent so applying $F$, we get in $X'$

\begin{picture}(30,28)(-9,-4)
\diagbw{}{}{}{}{}{}{}{} \diagbv{}{}{}{}{}{}{}{$\phi$}{}
\diagbr{}{}{}{}{$Ff$}{}{}{$\sim$}{$\phi$} \diagby{1mm}{$Fr$}
\end{picture}
\begin{picture}(12,18) \put(5,7){=} \end{picture}
\begin{picture}(25,15)(3,-2)
\onetwo{8}{0}{0.7mm}{}{}{$Ff$}{$Ff$}{$1$}
\end{picture}
\begin{picture}(12,18) \put(5,7){=} \end{picture}
\begin{picture}(30,22)(-9,-4)
\diagbw{}{}{}{}{$Ff$}{}{}{$\sim$} \diagby{1mm}{$r$}
\end{picture}

\begin{picture}(40,40)(-9,-14) \put(7,18){$\parallel$}
\diagbr{}{}{}{}{}{}{}{$\sim$}{$1 \ast \phi$} \diagby{1mm}{$\phi$}
\diagbz{$Fr$}
\end{picture}

\noindent as required.  The axiom for $l$ holds similarly.
\end{prf}

\addcontentsline{toc}{section}{References}
\bibliography{bib0209}

\nocite{bd2}
\nocite{hmp2} \nocite{hmp3} \nocite{hmp4}

\nocite{bae1}

\nocite{ben1}

\end{document}